\documentclass[conference]{IEEEtran}
\IEEEoverridecommandlockouts

\usepackage{amssymb}
\usepackage{algorithm}
\usepackage{algorithmicx}
\usepackage[noend]{algpseudocode}
\usepackage{bm}
\usepackage{booktabs}
\usepackage{capt-of}
\usepackage{color, colortbl}
\usepackage{mathtools}
\usepackage{multirow}
\usepackage{siunitx}
\usepackage{subcaption}
\usepackage[sort,compress]{cite}
\usepackage{xcolor}

\usepackage{lipsum}
\usepackage{amsfonts}
\usepackage{graphicx}
\usepackage{epstopdf}
\usepackage{hyperref}
\usepackage{cleveref}
\usepackage{listings}

\lstdefinelanguage{Julia}%
  {morekeywords={abstract,break,case,catch,const,continue,do,else,elseif,%
      end,export,false,for,function,immutable,import,importall,if,in,%
      macro,module,otherwise,quote,return,switch,true,try,type,typealias,%
      using,while},%
   sensitive=true,%
   alsoother={$},%
   morecomment=[l]\#,%
   morecomment=[n]{\#=}{=\#},%
   morestring=[s]{"}{"},%
   morestring=[m]{'}{'},%
}[keywords,comments,strings]%

\ifpdf
  \DeclareGraphicsExtensions{.eps,.pdf,.png,.jpg}
\else
  \DeclareGraphicsExtensions{.eps}
\fi

\newtheorem{remark}{Remark}

\newtheorem{assumption}{Assumption}
\newtheorem{definition}{Definition}
\newtheorem{theorem}{Theorem}

\Crefname{equation}{Eq.}{Eqns.}

\title{
{An O(N) distributed-memory parallel direct solver for planar integral equations}
}

\makeatletter
\newcommand{\linebreakand}{%
  \end{@IEEEauthorhalign}
  \hfill\mbox{}\par
  \mbox{}\hfill\begin{@IEEEauthorhalign}
}
\makeatother

\author{\IEEEauthorblockN{Tianyu Liang}
\IEEEauthorblockA{University of California, Berkeley, USA \\
tianyul@berkeley.edu}
\and
\IEEEauthorblockN{Chao Chen}
\IEEEauthorblockA{North Carolina State University, USA \\
cchen49@ncsu.edu} \\
\linebreakand
\IEEEauthorblockN{Per-Gunnar Martinsson}
\IEEEauthorblockA{The University of Texas at Austin, USA \\
pgm@oden.utexas.edu}
\and
\IEEEauthorblockN{George Biros}
\IEEEauthorblockA{The University of Texas at Austin, USA \\
biros@oden.utexas.edu}
}

\usepackage{amsopn}

\usepackage{booktabs}
\usepackage{multirow}
\usepackage{siunitx}

\ifpdf
\hypersetup{
  pdftitle={Distributed IE solver},
  pdfauthor={T. Liang, C. Chen, G. Martinsson, and G. Biros}
}
\fi

\graphicspath{{figs/}}

\newcommand{\N}{ {\mathcal{N}} }

\newcommand{\bigO}{ {\mathcal{O}} }
\newcommand{\B}{ {\mathcal{B}} }
\newcommand{\X}{ {\mathcal{X}} }
\newcommand{\F}{ {\mathcal{F}} }
\newcommand{\T}{ {\mathcal{T}} }
\newcommand{\eps}{ {\varepsilon} }

\newcommand{\J}{ {\mathcal{J}} }
\newcommand{\R}{ {\mathcal{R}} }
\newcommand{\C}{ {\mathcal{C}} }
\newcommand{\M}{ {\mathcal{M}} }

\renewcommand{\S}{ {\mathcal{S}} }
\newcommand{\W}{ {\mathcal{W}} }

\newcommand{\cyan}[1]{\textcolor{cyan}{#1}}

\newcolumntype{H}{>{\setbox0=\hbox\bgroup}c<{\egroup}@{}}

\definecolor{Gray}{gray}{0.9}

\begin{document}

\maketitle

\begin{abstract}
Boundary value problems involving elliptic PDEs such as the Laplace and the Helmholtz equations are ubiquitous in mathematical physics and engineering. Many such problems can be alternatively formulated as integral equations that are mathematically more tractable. However, an integral-equation formulation poses a significant computational challenge: solving large dense linear systems that arise upon discretization. In cases where iterative methods converge rapidly, existing methods that draw on fast summation schemes such as the Fast Multipole Method are highly efficient and well-established. More recently, linear complexity direct solvers that sidestep convergence issues by directly computing an invertible factorization have been developed. However, storage and computation costs are high, which limits their ability to solve large-scale problems in practice. In this work, we introduce a distributed-memory parallel algorithm based on an existing direct solver named ``strong recursive skeletonization factorization~\cite{minden2017recursive}.'' 
Specifically, we apply low-rank compression to certain off-diagonal matrix blocks in a way that minimizes computation and data movement. 
Compared to iterative algorithms, our method is particularly suitable for problems involving ill-conditioned matrices or multiple right-hand sides. Large-scale numerical experiments are presented to show the performance of our \texttt{Julia}  implementation.

\end{abstract}



\section{Introduction}

Boundary value problems of classical potential theory appear frequently in scientific and engineering domains. While it may be challenging to solve these problems directly through their typical formulations involving elliptic partial differential equations (e.g.,  Laplace's and the Helmholtz equations), alternative formulations of these problems via integral equations (IEs) can be solved more efficiently in many environments~\cite[Chapter~10]{martinsson2019fast}. An integral equation takes the  form
\begin{equation} \label{e:ie}
a(\bm{x}) u(\bm{x}) + \int_{\Omega} K(\bm{x}, \bm{y}) u(\bm{y}) d \bm{y} = f(\bm{x}), \quad \bm{x} \in \Omega
\end{equation}
where $u(\bm{x})$ is the unknown to be determined, $K(\bm{x}, \bm{y})$ is derived from the free-space fundamental solution associated with the underlying elliptic operator, $a(\bm{x})$ and $f(\bm{x})$ are given functions,
and  $\Omega$ is the problem domain. 

We are interested in solving the associated \emph{dense} linear system
\begin{equation} \label{e:axb}
A \, x = b, \quad A \in \mathbb{C}^{N \times N}, \quad x \, \text{ and } \, b  \in \mathbb{C}^{N}
\end{equation}
that arises from the discretization of \Cref{e:ie} using  approaches such as collocation, the Nystr\"om method, or the Galerkin method.
To construct a fast solver, we exploit the fact that $A$ is a kernel matrix, that is
\begin{equation} \label{e:aij}
A_{i,j} = K (\bm{x}_i, \bm{x}_j), \quad \forall i \not = j,
\end{equation}
where $\{\bm{x}_i\}_{i=1}^N \subset \Omega$ are points related to the discretization of \Cref{e:ie}. Here, the kernel function $K$ comes from \Cref{e:ie} and is smooth except when $\bm{x}_i = \bm{x}_j$. In this paper, we restrict ourselves to the case that $K$ is not too oscillatory and the problem domain $\Omega$ is planar.

Since  kernel function $K$  is smooth in \Cref{e:aij}, certain off-diagonal blocks in $A$ can be compressed using low-rank approximations. In particular, the required numerical rank is constant (independent of $N$) for a prescribed accuracy when an off-diagonal block satisfies the so-called \emph{strong admissibility}. This observation is crucial to the success of the fast multipole method (FMM)~\cite{greengard1987fast,greengard1996new}, which requires only $\bigO(N)$ operations for applying matrix $A$ to a vector. It may be somewhat surprising that under mild assumptions motivated by numerical experiments, a factorization of $A$ can also be computed using $\bigO(N)$ operations~\cite{ambikasaran2014inverse,corona2015n,ho2016hierarchical,minden2017recursive,coulier2017inverse}. Although the asymptotic complexity is appealing, such factorization algorithms still require a significant amount of computation and memory footprint. A natural solution is to distribute the required computation and storage across multiple compute nodes in order to solve large-scale problems. The challenge here comes from the fact that  the matrix $A$ is dense, which seems to imply that the cost of communication among compute nodes may be excessive. 



\subsection{Contributions}

{We introduce a distributed-memory parallel {direct solver} based on the  (sequential) strong recursive skeletonization (RS-S)~\cite{minden2017recursive}, {one of a few methods that require $\bigO(N)$ storage and computation} for solving \emph{dense} linear systems \Cref{e:axb} arising from the discretization of planar integral equations. 
To distribute the storage and the computation, we partition the problem domain $\Omega$ among a grid of $p$ processes, and we show that in our distributed-memory algorithm every process needs to communicate with only its neighbors in the process grid. As a result, every process sends a total of $\bigO(\log N + \log p)$ messages with a total of $\bigO(\sqrt{N/p} + \log p)$ words (under the assumptions in \Cref{s:analysis}).
We also describe an implementation using  the \texttt{Julia} programming language~\cite{bezanson2012julia,bezanson2017julia} that can tackle dense matrices with 1 \emph{billion} rows and columns.}

{As a direct solver,} our algorithm has two stages: (1) a factorization stage that constructs a compressed form of the inverse $A^{-1}$ subject to a prescribed accuracy and (2) a solution stage that applies the compressed inverse to a right-hand side  (vector) $b$ to obtain an approximate solution of \Cref{e:axb}. 
Once the first phase finishes, the second phase is extremely efficient. Therefore, our algorithm is particularly suitable for problems where
\begin{enumerate}
\item
the condition number of  $A$ is large. This typically occurs when \Cref{e:ie} is a first-kind Fredholm integral equation for Laplace’s equation or the Stokes equation, which has applications in magnetostatics, electrostatics and fluid dynamics. For example, the condition number of $A$ scales empirically as $\bigO(N)$ for Laplace’s equation.

\item 
multiple right-hand sides need to be solved. 
This typically occurs when  \Cref{e:ie}  is the Lippmann-Schwinger equation that models the scattering of acoustic waves in media with variable speed. Solving the Lippmann-Schwinger equation has applications in seismology, ultrasound tomography, and sonar, where incident waves from multiple angles need to be solved. 
 
\end{enumerate}


As in the FMM, our method starts by partitioning the set of row/column indices $\{1,2,\ldots,N\}$ of matrix $A$  into clusters $\B_1, \B_2, \ldots$ such that indices in every cluster correspond to points spatially close to each other. The mathematical property we exploit in our method is the following.
Given a cluster $\B$ and the union of non-adjacent clusters $\F$, the off-diagonal submatrices $A_{\mathcal B, \mathcal F}$ and $A_{\mathcal F, \mathcal B}$ are numerically low-rank and can be compressed efficiently without explicitly forming the entire matrices.
We show that applying Gaussian elimination to a subset of indices in $\B$ after the compression (a.k.a, skeletonization) leads to a Schur-complement update that affects only  adjacent clusters.

Based on these observations, we propose a domain decomposition strategy where the computational domain $\Omega$ is partitioned among all processes and show that all processes can work on their interior clusters in parallel. To process the remaining boundary clusters, we  color all processes so any pair of  processes have  two different colors if they own adjacent clusters. Then, we loop over all colors, and processes with the same color can work in parallel at each iteration.

We implemented our parallel algorithm using \texttt{Julia}  because of its ease of usage and support for various numerical linear algebra routines. In particular, our implementation explores the distributed computing capability of \texttt{Julia}\footnote{\url{https://docs.julialang.org/en/v1/manual/distributed-computing/}}. 



\subsection{{Relation to existing work}}

{The only existing distributed-memory parallel direct solver with $\bigO(N)$ complexity for solving \Cref{e:axb}  is  Ma et al.~\cite{Ma2022}.} They introduced an $\mathcal{H}^2$-ULV factorization in two passes: the first pass pre-computes all possible fill-ins and includes them in  low-rank approximations, and thus the second pass computing the actual factorization is fully parallel. {The disadvantage is that the overhead of an extra pass can be significant.}

The inverse FMM~\cite{ambikasaran2014inverse,coulier2017inverse} is another   $\bigO(N)$ (sequential) method and shares a lot of similarities with the RS-S~\cite{minden2017recursive}, which our distributed-memory algorithm is built upon. Takahashi et al.~\cite{takahashi2020parallelization} analyzed the parallelism of the inverse FMM and proposed a parallel algorithm for \emph{shared-memory} machines. {The proposed algorithm by Takahashi et al.~\cite{takahashi2020parallelization} is similar to a parallel strategy briefly discussed at the end of~\cite{minden2017recursive}, which colors all clusters of indices (our distributed-memory algorithm colors processes).}

Both the inverse FMM and the RS-S implicitly make use of the {$\mathcal{H}^2$-matrices introduced by Hackbusch and collaborators.} Their seminal work~\cite{hackbusch1999sparse,hackbusch2002data,5555} establishes the algebra of $\mathcal{H}$- and $\mathcal{H}^2$-matrices, which was implemented by several software packages~\cite{grasedyck2008parallel,kriemann2013h,kriemann2005parallel} attaining the theoretical linear or quasilinear complexity for solving \Cref{e:axb}. {The  algorithms are recursive and rely on expensive hierarchical matrix-matrix multiplication, so the hidden constants in the asymptotic scalings are quite large.}

Other forms of hierarchical low-rank approximations include the hierarchical semi-separable (HSS) matrices~\cite{chandrasekaran2006fast,chandrasekaran2007fast,xia2010superfast}, hierarchical off-diagonal low-rank (HODLR) matrices~\cite{ambikasaran2013mathcal,aminfar2016fast}, among others~\cite{chen2002fast,bremer2012fast}. These methods have $\bigO(N)$ complexity in 1D (e.g., boundary IEs on curves), {but their costs deteriorate to  super linear in 2D and 3D if the a fixed target accuracy is desired.} High-performance and parallel implementations can be found in~\cite{rouet2016distributed,ghysels2016efficient,liu2016parallel,cai2018smash,ambikasaran2019hodlrlib,chenhan2019distributed,chen2022solving,deshmukh2023n}.
Corona et al.~\cite{corona2015n} and Ho et al.~\cite{ho2016hierarchical}  improve upon earlier recursive skeletonization ideas~\cite{martinsson2005fast,greengard2009fast} and attain $\bigO(N)$ complexity by incorporating  extra compression steps on intermediate skeletons. However,  parallelizing such algorithms is challenging.

Besides approximations in hierarchical formats,  {flat formats such as block low-rank~\cite{mary2017block} or tile low-rank~\cite{akbudak2017tile}} have also been proposed recently. Although algorithms leveraging these flat formats have led to significant speedups compared to classical algorithms for general dense matrices in practical applications~\cite{amestoy2019performance,boukaram2021h2opus,cao2022framework,Cao2022}, {they do not attain $\bigO(N)$ complexity.}

\section{Sequential algorithm and data dependency}

In this section, we review the sequential algorithm for factorizing and solving \Cref{e:axb} approximately. The focus here is analyzing data dependency and parallelism. While there exist a few  methods that are closely related such as~\cite{ambikasaran2014inverse,coulier2017inverse,pouransari2017fast,sushnikova2018compress}, we mostly follow~\cite{minden2017recursive,sushnikova2022fmm}.

Suppose we are given a set of points $\mathcal{X} = \{\bm{x}_i\}_{i=1}^N$ in $\mathbb{R}^{2}$ and a kernel function $K$ that defines matrix $A$ as in \Cref{e:aij}. Our approach for solving \Cref{e:axb} employs approximate Gaussian elimination in a multi-level fashion.

\subsection{Hierarchical domain decomposition} \label{s:dd}

\begin{figure}
\centering
     \begin{subfigure}{0.23\textwidth}
         \centering
         \includegraphics[width=\linewidth]{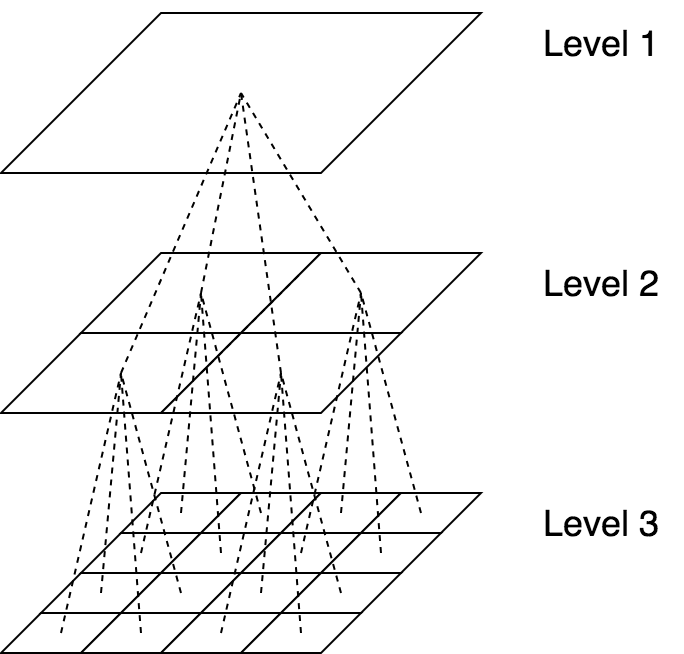}
         \caption{Hierarchical domain decomposition and the quad-tree.}
         \label{f:tree}
     \end{subfigure}
     \hfill
     \begin{subfigure}{0.23\textwidth}
         \centering         
         \includegraphics[trim=0 -30 0 0,clip,width=\linewidth]{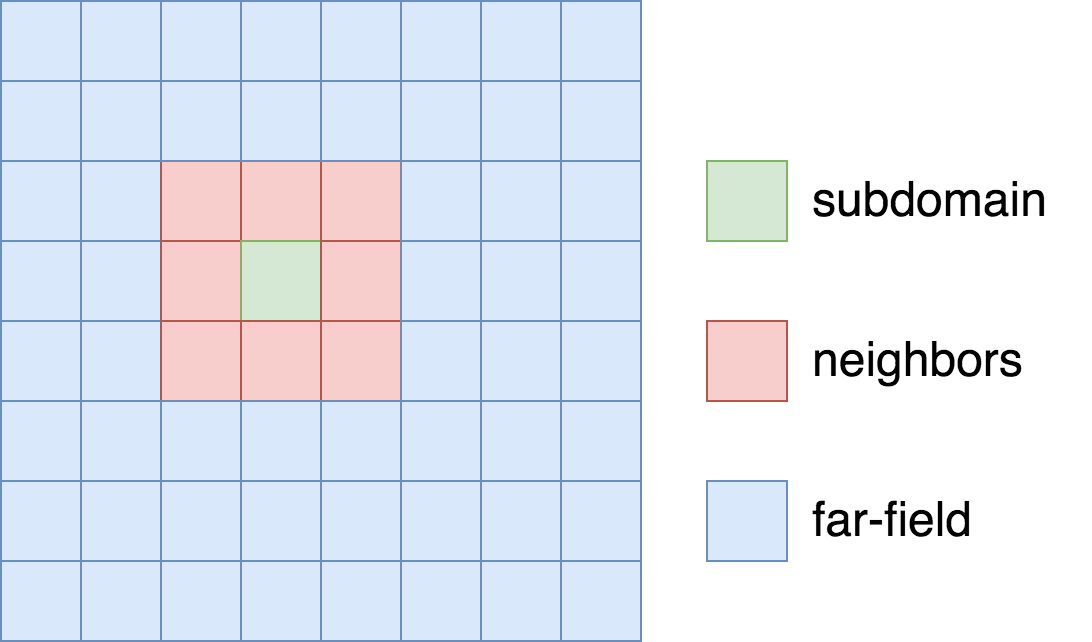}
         \caption{A subdomain, its near-field (neighbors), and its far-field.}
         \label{f:bnf}
     \end{subfigure}
     \caption{Assume the points are uniformly distributed in a square domain. (a) a 3-level hierarchical domain decomposition and the corresponding quad-tree. (b) a subdomain at the fourth  level in a hierarchical domain decomposition, its near-field (neighbors), and its far-field.}
\end{figure}

The algorithm relies on a hierarchical decomposition of the points $\mathcal{X} = \{\bm{x}_i\}_{i=1}^N$ as in the FMM. The hierarchical decomposition is as follows. Suppose all points lie in a square domain. We divide the domain into four equal-sized subdomains. We call the entire domain the \emph{parent} of the four subdomains, which are the \emph{children} of the original domain. We continue subdividing the four subdomains recursively until every subdomain contains $\bigO(1)$ points. This hierarchical decomposition of the problem domain naturally maps to a  quad tree $\mathcal{T}$, where the root  is the entire domain and the leaves are all subdomains that are not subdivided. See a pictorial illustration  in  \Cref{f:tree}.

For ease of presentation, we further make the assumption that the points are uniformly distributed, so the hierarchical domain decomposition corresponds to a \emph{perfect} quad tree (all internal nodes have 4 children and all the leaf nodes are at the same level). Extensions to a non-uniform distribution of points are straightforward but quite tedious; see details in~\cite{minden2017recursive,sushnikova2022fmm}.

We refer to a  subdomain (or a  node in tree $\mathcal{T}$) as a \emph{box} for the rest of this paper. Given a box $\mathcal{B}$ at level $\ell$ in the tree $\mathcal{T}$, its \emph{parent} $\mathcal{P}(\B)$ and its  \emph{children} $\mathcal{C}(\B)$ are naturally defined (if they exist). We also define the \emph{near-field} (a.k.a., \emph{neighbors}) $\mathcal{N}(\B)$ and the \emph{far-field} $\mathcal{F}(\B)$ as follows:
\begin{itemize}
%
\item
$\mathcal{N}(\B)$: boxes that are physically adjacent to $\B$ (excluding $\B$) at level $\ell$.

\item
$\mathcal{F}(\B)$:  boxes excluding $\B \cup \mathcal{N}(\B)$ at level $\ell$.

\end{itemize}
See an example  in  \Cref{f:bnf}. Assuming $\T$ is a perfect quad-tree, no box has  more than eight  neighbors, i.e., $|\N(\B)| \le 8$.
With some  abuse in notation, we use $\B, \N, \F$ to also denote indices of points in $\X$ that are contained in box $\B$, its near-field/neighbors $\N(\B)$, and its far-field $\F(\B)$, respectively.

Given a leaf node $\mathcal B$, there exists a permutation matrix $P$ such that
\begin{equation}
    \label{ty_step1}
    P^{\top}AP = \begin{pmatrix}
    A_{\mathcal B, \mathcal B} & A_{\mathcal B, \mathcal{N}} & A_{\mathcal B, \mathcal F} \\
    A_{\mathcal{N}, \mathcal B} & A_{\mathcal{N}, \mathcal{N}} & A_{\mathcal{N}, \mathcal F} \\
    A_{\mathcal F, \mathcal B} & A_{\mathcal F, \mathcal{N}} & A_{\mathcal F, \mathcal F}
    \end{pmatrix}.
\end{equation}
Here, $A_{\mathcal I, \mathcal J}$ for two index sets $\mathcal I$ and $\mathcal J$ stands for the corresponding submatrix in $A$.

\subsection{Low-rank property} \label{s:lr}

\begin{assumption} \label{a:bf}
For any box $\B$, it holds that 
$
A_{i,j} = K(x_i, x_j) 
$
for (1) $i \in \B$ and $j \in\F$; or (2) $i \in \F$ and $j \in\B$.
In other words, the submatrices $A_{\mathcal B, \mathcal F}$ and $A_{\mathcal F, \mathcal B}$ come from kernel evaluation.
\end{assumption}
%

The  observation is that the off-diagonal block $A_{\mathcal B, \mathcal F}$ (or $A_{\mathcal F, \mathcal B}$) for an arbitrary box $\B$ can be  approximated efficiently by a low-rank approximation for a prescribed accuracy $\eps$. 
In particular, the numerical rank is $\bigO(1)$, independent of the problem size $N$, for a smooth kernel that is not highly oscillatory.


Next, we explain a specific type of low-rank approximation named the interpolative decomposition (ID)~\cite{cheng2005compression}, which is applied to the submatrix $A_{\mathcal B, \mathcal F}$. Other low-rank approximations such as the truncated  singular value decomposition can also be used; see the resulting solvers for \Cref{e:axb} in, e.g.,~\cite{pouransari2017fast,cambier2020algebraic}.

\begin{definition}
\label{definition:id}
Let  $\mathcal{I} = \{1,2,\ldots,m\}$ and $\mathcal{J} = \{1,2,\ldots,n\}$ be the row and column indices of a matrix $A_{\mathcal I, \mathcal J} \in \mathbb{C}^{m \times n}$. A (column) interpolative decomposition (ID) for a prescribed accuracy $\eps$ finds the so-called \emph{skeleton} indices $\S \subset \J$, the \emph{redundant} indices $\R = \J \backslash \S$, and an \emph{interpolation matrix} $T \in \mathbb{C}^{|\S| \times |\R|}$ such that
\[
    \|A_{\mathcal I, \mathcal R} - A_{\mathcal I, \mathcal S} \, T\| 
     \le \eps \, \|A_{\mathcal I, \mathcal J} \|.
\]
\end{definition}

While the strong rank-revealing QR factorization of Gu and Eisenstat~\cite{gu1996efficient} is the most robust method for computing an ID, we employ   the column-pivoting QR factorization as a greedy approach~\cite{cheng2005compression}, which has better computational efficiency and  behave well in practice. In particular, we adopted the implementation from the \texttt{LowRankApprox.jl}\footnote{\url{https://github.com/JuliaLinearAlgebra/LowRankApprox.jl}} package. The  cost to compute an ID using the aforementioned deterministic methods is $\bigO(mn \,|\S|)$, which can be further reduced to $\bigO(mn \,\log(|\S|) + |\S|^2 n)$ using randomized algorithms that may incur some loss of accuracy~\cite{dong2021simpler}.

\subsection{Fast compression} \label{s:comp}

As stated earlier, we want to compress the two off-diagonal blocks $A_{\mathcal B, \mathcal F}$ and $A_{\mathcal F, \mathcal B}$ using their IDs. In practice, we conduct one (column) ID compression of the concatenation 
\begin{equation} \label{e:con}
\begin{pmatrix}
A_{\mathcal F, \mathcal B} \\
A_{\mathcal B, \mathcal F}^*
\end{pmatrix}
\end{equation}
(that has $\bigO(N)$ rows and $\bigO(1)$ columns) and obtain a single interpolation matrix $T$ that satisfies both
\begin{equation} \label{e:t}
A_{\mathcal F, \mathcal R} \approx A_{\mathcal F, \mathcal S} \, T
\quad \text{and} \quad
A_{\mathcal R, \mathcal F} \approx T^* \, A_{\mathcal S, \mathcal F} 
.
\end{equation}
This approach leads to a slightly larger set of skeleton indices but makes the algorithm/implementation easier.

Notice that \emph{the computational cost would be $\bigO(N)$ if the full matrix in \Cref{e:con} is formed, which turns out to be unnecessary}. As in the FMM, there are a few techniques that require only $\bigO(1)$ operations such as the (analytical) multipole expansion (see, e.g.,~\cite{greengard1987fast,greengard1996new,cheng2006wideband}), the Chebyshev interpolation (see, e.g.,~\cite{fong2009black,wang2021pbbfmm3d}), and the proxy method (see, e.g.,~\cite{ying2004kernel,martinsson2005fast}).

\begin{figure}
     \centering
     \begin{subfigure}{0.2\textwidth}
         \centering
         \includegraphics[width=\linewidth]{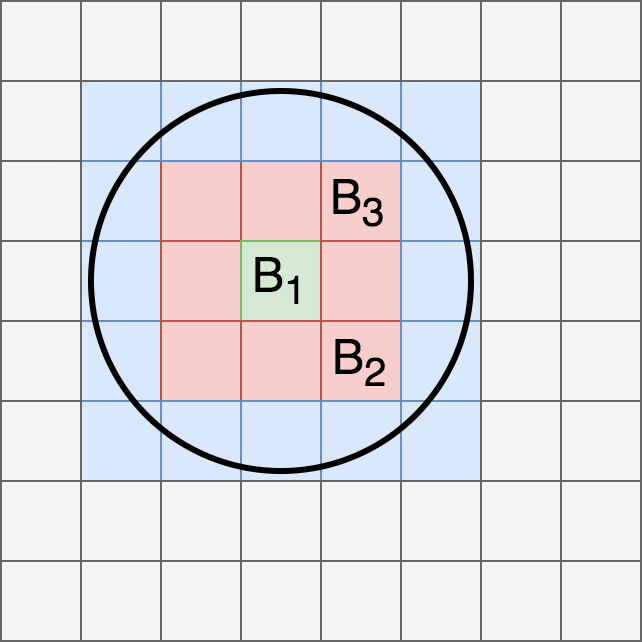}
         \caption{}
     \end{subfigure}
     \hfill
     \begin{subfigure}{0.2\textwidth}
         \centering         
         \includegraphics[width=\linewidth]{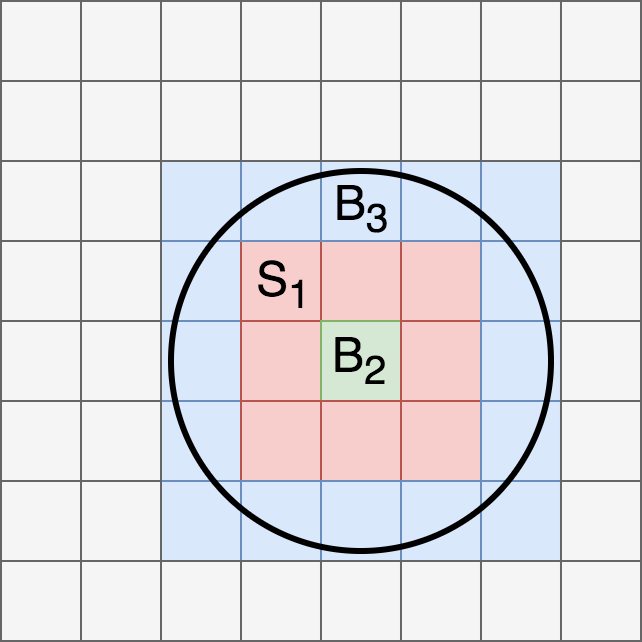}
         \caption{}
     \end{subfigure}
     \caption{Illustration of proxy circle and distance-2 neighbors $\M$. 
     (a) For a box $\B_1$, its distance-2 neighbors $\M(\B_1)$ (in blue) is the {small subset (16 boxes at most)}  of its far-field $\F(\B_1)$ that surrounds its neighbors (in red). The proxy circle for $\B_1$ lies in $\M(\B)$. 
     %
     (b) An example showing that the ``interactions''  between $\B_3$ and $\B_2$ (submatrices $A_{\B_2, \B_3}$ and $A_{\B_3, \B_2}$) have been modified after the redundant indidces $\R_1$ is eliminated (the skeleton indices $\S_1$ remain). See \Cref{r:nn} below.
     }
     \label{f:proxy}
\end{figure}

Below, we focus on illustrating the proxy method for compressing the submatrix $A_{\mathcal F, \mathcal B}$. Specifically, we form and compress the following matrix
\begin{equation} \label{e:proxy}
\begin{pmatrix}
A_{\M, \B} \\
K_{\text{proxy}, \B}
\end{pmatrix}
,
\end{equation}
which has \emph{$\bigO(1)$ rows} (and columns). In \Cref{e:proxy}, $\M = \M(\B)$ is a \emph{small} subset of $\F(\B)$ surrounding $\N(\B)$  defined as follows (see \Cref{f:proxy}):

\begin{definition} \label{d:m}
The distance-2 neighbors of a box $\B$ is defined as
\[
\M(\B) = \N( \N(\B) ) \backslash \left(\N(\B) \cup \B \right).
\]
\end{definition}

Note that \Cref{a:bf} does not  hold in our algorithm, where the submatrix $A_{\M, \B}$ is updated and is \emph{not} the kernel evaluation any more; see \Cref{s:nn}. However, we can show the following  to be true:

\begin{theorem} \label{a:bm}
For any leaf box $\B$, it holds that 
$
A_{i,j} = K(x_i, x_j) 
$
for (1) $i \in \B$ and $j \in\F\backslash \M$; or (2) $i \in \F\backslash \M$ and $j \in\B$.
In other words, the submatrices $A_{\mathcal B, \mathcal F \backslash \M}$ and $A_{\mathcal F \backslash \M, \mathcal B}$ come from kernel evaluation.
\end{theorem}

This relaxation means that the established rank estimation in \Cref{s:lr} no longer holds. However, we empirically verify that the resulting rank from compressing \Cref{e:proxy} still follows the previous estimation; see numerical results in \Cref{s:r}.

In \Cref{e:proxy}, $K_{\text{proxy}, \B}$ stands for a matrix where the $(i,j)$-th entry is given by $K(y_i, x_j)$, i.e., a kernel evaluation between a discretization point $y_i$ on the \emph{proxy circle}\footnote{A circle is chosen for ease of implementation but not necessary mathematically.} and a point $x_j \in \B$ lying inside the box $\B$. 
The ``interaction'' $K_{\text{proxy}, \B}$ accounts for the  ``interaction'' between points lying in $\F(\B) \backslash \M(\B)$ and points inside $\B$. To that end, the proxy circle must lie inside $ \M(\B)$. See a pictorial illustration in \Cref{f:proxy}. In this paper, we choose the radius of the  proxy circle to be $2.5 L$, where $L$ is the side length of boxes at this level.

\begin{remark} \label{r:t}
The computation of skeleton indices $\S$, redundant indices $\R$, and the interpolation matrix $T$ in \Cref{e:t} requires (reads) submatrices $A_{\M, \B}$ and $A_{\B, \M}$ (\emph{not} the entire submatrices $A_{\F, \B}$ or $A_{\B, \F}$) for a leaf box $\B$.
\end{remark}

\subsection{Approximate matrix factorization} \label{s:nn}

The algorithm starts at the leaf level of the quad-tree $\mathcal T$. Let us apply an ID compression for a box $\B$ to obtain skeleton indices $\S$, redundant indices $\R = \B \backslash \S$, and an interpolation matrix $T$ such that \Cref{e:t} holds. (As explained in the previous section, this step requires only $\bigO(1)$ operations.)

Rewrite  \Cref{ty_step1} as 
\begin{align*}
    \label{ty_low_rank}
     \tilde{P}^\top A \tilde{P} &=
     \left(
     \begin{array}{cc|cc}
    A_{\mathcal R, \mathcal R} & A_{\mathcal R, \mathcal S} & A_{\mathcal R, \mathcal{N}} & A_{\mathcal R, \mathcal F} \\
    A_{\mathcal S, \mathcal R} & A_{\mathcal S, \mathcal S} & A_{\mathcal S, \mathcal N} & A_{\mathcal S, \mathcal F} \\
    \hline
    A_{\mathcal N, \mathcal R} & A_{\mathcal N, \mathcal S} & A_{\mathcal N, \mathcal N} & A_{\mathcal N, \mathcal F} \\
    A_{\mathcal F, \mathcal R} & A_{\mathcal F, \mathcal S} & A_{\mathcal F, \mathcal N} & A_{\mathcal F, \mathcal F}\\
    \end{array}
    \right) \\
        & 
        \stackrel{\Cref{e:t}}{\approx}
     \left(
     \begin{array}{cc|cc}
    A_{\mathcal R, \mathcal R} & A_{\mathcal R, \mathcal S} & A_{\mathcal R, \mathcal{N}} & T^{*}A_{\mathcal S, \mathcal F} \\
    A_{\mathcal S, \mathcal R} & A_{\mathcal S, \mathcal S} & A_{\mathcal S, \mathcal N} & A_{\mathcal S, \mathcal F} \\
    \hline
    A_{\mathcal N, \mathcal R} & A_{\mathcal N, \mathcal S} & A_{\mathcal N, \mathcal N} & A_{\mathcal N, \mathcal F} \\
    A_{\mathcal F, \mathcal S} \, T & A_{\mathcal F, \mathcal S} & A_{\mathcal F, \mathcal N} & A_{\mathcal F, \mathcal F}\\
    \end{array}
    \right),
\end{align*}
where $\tilde{P}$ is an appropriate permutation matrix. Define the \emph{sparsification} matrix
\begin{equation} \label{e:s}
S = 
\begin{pmatrix}
I \\
-T & I \\
&& I \\
&&& I
\end{pmatrix},
\end{equation}
where the partitioning of row and column indices is the same as that of $\tilde{P}^\top A \tilde{P} $ and the diagonal blocks are identity matrices of appropriate sizes. Applying the sparsification matrix, we have that
\[
S^* \, (\tilde{P}^\top A \tilde{P}) \, S  \approx
     \left(
     \begin{array}{cc|cc}
    X_{\mathcal R, \mathcal R} & X_{\mathcal R, \mathcal S} & X_{\mathcal R, \mathcal{N}}  
    \\
    X_{\mathcal S, \mathcal R} & A_{\mathcal S, \mathcal S} & A_{\mathcal S, \mathcal N} & A_{\mathcal S, \mathcal F} \\
    \hline
    X_{\mathcal N, \mathcal R} & A_{\mathcal N, \mathcal S} & A_{\mathcal N, \mathcal N} & A_{\mathcal N, \mathcal F} \\
     & A_{\mathcal F, \mathcal S} & A_{\mathcal F, \mathcal N} & A_{\mathcal F, \mathcal F}\\
    \end{array}
    \right),
\]
where {the coupling (submatrices) between $\R$ and $\F$ disappears}. Here, the notation $X$ denotes modified blocks that can be easily derived and be calculated using $\bigO(1)$ operations. \emph{Notice that we never need to explicitly form $A_{\mathcal R, \mathcal F}$ or $A_{\mathcal S, \mathcal F}$.}

Suppose that the diagonal block $X_{\R, \R}$ is non-singular and that 
$X_{\R, \R} = L_{\R} U_{\R}$
is an LU factorization. We apply (block) Gaussian elimination to obtain
\begin{equation}
    \label{e:step}
     L \, ( S^* \tilde{P}^\top A \tilde{P} S ) \, U \approx
     \left(
     \begin{array}{cc|cc}
    I &
    \\
     & X_{\mathcal S, \mathcal S} & X_{\mathcal S, \mathcal N} & A_{\mathcal S, \mathcal F} \\
    \hline
     & X_{\mathcal N, \mathcal S} & X_{\mathcal N, \mathcal N} & A_{\mathcal N, \mathcal F} \\
     & A_{\mathcal F, \mathcal S} & A_{\mathcal F, \mathcal N} & A_{\mathcal F, \mathcal F}\\
    \end{array}
    \right),
\end{equation}
where 
\begin{align*}
& L = 
\begin{pmatrix}
I \\
-X_{\S, \R} U_{\R}^{-1} & I \\
-X_{\N, \R} U_{\R}^{-1} && I \\
&&& I
\end{pmatrix}
\begin{pmatrix}
L_{\R}^{-1} \\
& I \\
&& I \\
&&& I
\end{pmatrix}
\quad \text{and} \quad  \\
& U = 
\begin{pmatrix}
U_{\R}^{-1} \\
& I \\
&& I \\
&&& I
\end{pmatrix}
\begin{pmatrix}
I & - L_{\R}^{-1}  X_{\R, \S} & - L_{\R}^{-1} X_{\R, \N}  \\
& I \\
 && I \\
&&& I
\end{pmatrix}.
\end{align*}
Notice that \emph{the original ``interaction'' among $\B$'s neighbors, i.e., submatrix $A_{\mathcal N, \mathcal N}$, has been updated.}

Define 
\begin{equation} \label{e:sk}
V = L S^* \tilde{P}^\top,
\quad \quad
W = \tilde{P} S U,
\quad \text{and} \quad
Z(A; \B) = V A W.
\end{equation}
Here, $Z(A; \B)$ is called the \emph{strong skeletonization operator} in~\cite{minden2017recursive}. Note that $V^{-1}$ and $W^{-1}$ can be applied efficiently without explicitly forming their inverses.

\begin{remark} \label{r:nn}
Given the interpolation matrix $T$ in \Cref{e:t}, applying the strong skeletonization operator $Z(A; \B)$ requires (reads) matrices $X_{\N, \R}$ and $X_{\R, \N}$, and it updates (read \& write) submatrix $A_{\mathcal N, \mathcal N}$ for a leaf box $\B$. (No far-field information is required.)
\end{remark}

\subsection{Multi-level algorithm} \label{s:ml}

Let $n_\ell$ denote the number of boxes at level $\ell$ in the quad-tree $\T$ ($\ell=1,2,\ldots,L$). Given an ordering of all boxes $\B_1,\B_2,\ldots, \B_{n_L}$ at the leaf level, our algorithm applies the strong skeletonization operator \Cref{e:sk} to all boxes one after another:
\begin{align*}
Z(A; \B_1, \B_2) & \triangleq Z( Z(A; \B_1); \B_2) \\
Z(A; \B_1, \B_2, \B_3) & \triangleq Z( Z(A; \B_1, \B_2); \B_3) \\
& \ldots
\end{align*}
The resulting (approximate) factorization becomes
\begin{equation} \label{e:level}
Z(A; \B_1,\B_2,\ldots, \B_{n_L}) \approx 
\begin{pmatrix}
I_{\R_1} \\
& I_{\R_2} \\
&& \ddots \\
&&& I_{\R_{n_L}} \\
&&&& \tilde{A} \\
\end{pmatrix}
,
\end{equation}
where the leading diagonal blocks are identity matrices of sizes $|\R_i|$, which correspond to the redundant indices in box $\B_i$ for $i=1,2,\ldots,n_L$. The (approximate) Schur complement $\tilde{A}$ has $\sum_i \S_i$ rows and columns.

\begin{figure}
     \centering
     \begin{subfigure}{0.25\textwidth}
         \centering
         \includegraphics[width=\linewidth]{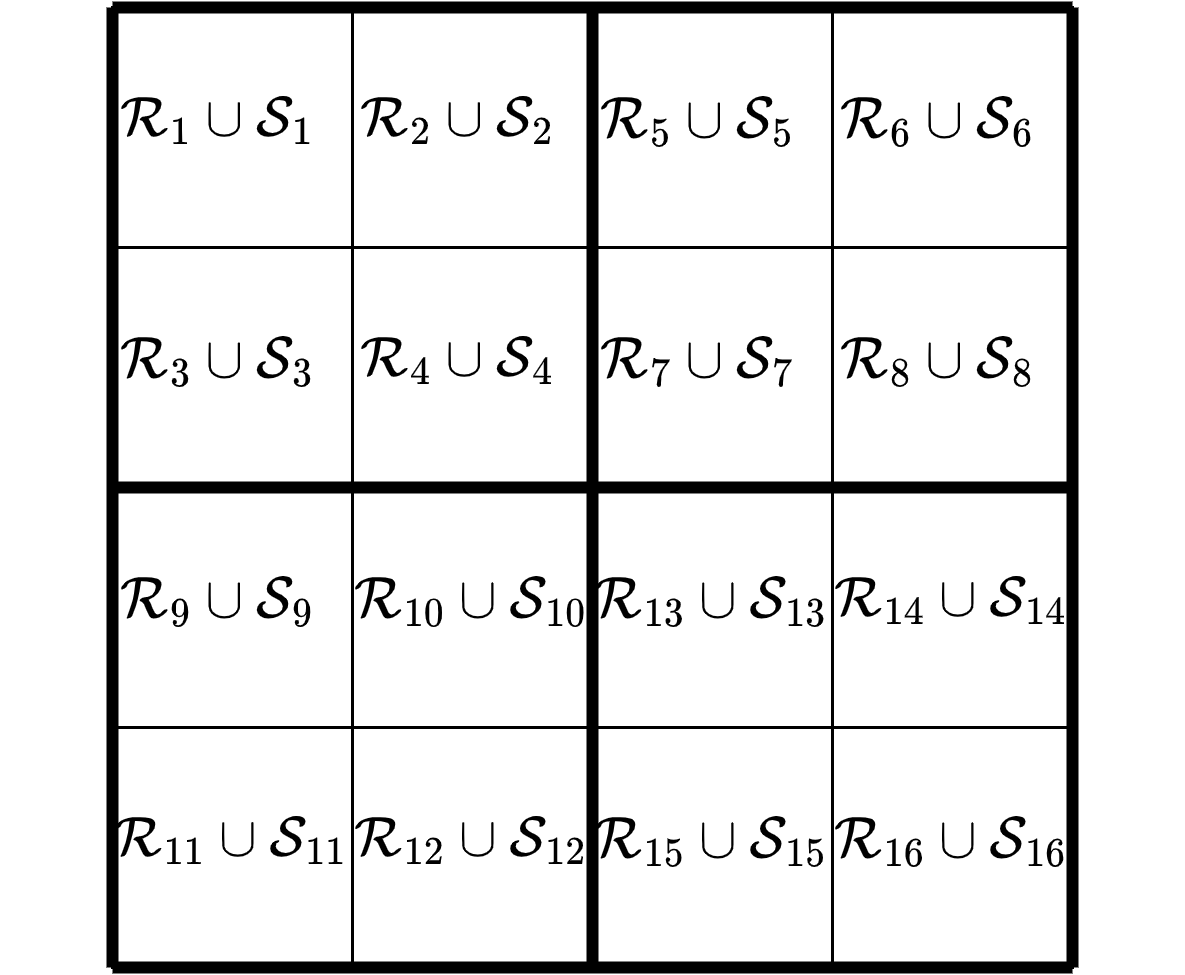}
         \caption{The leaf level.}
     \end{subfigure}
     \hfill
     \begin{subfigure}{0.214\textwidth}
         \centering         
         \includegraphics[width=\linewidth]{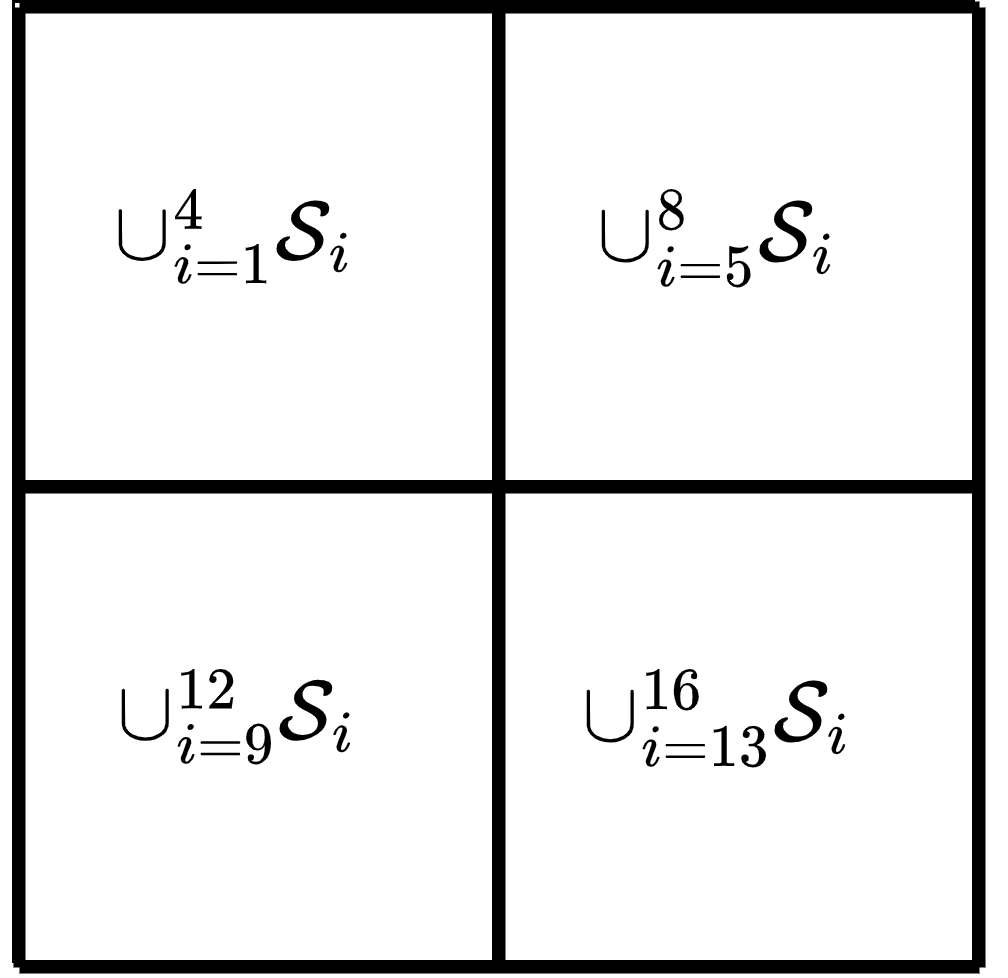}
         \caption{The parent level.}
     \end{subfigure}
     \caption{Illustration of the merge process. Suppose a square domain is divided into 16 boxes $\B_1, \B_2, \ldots, \B_{16}$. 
     (a) At the leaf level, a box $\B_i$ owns points/indices that physically lie inside it. These indices are divided into skeleton indices $\S_i$ and redundant indices $\R_i$ after an ID compression. (b) At the coarse/parent level, a parent box owns the skeleton indices of its children. 
     (Again, they will be divided into skeleton and redundant indices after an ID compression.)}
     \label{f:merge}
\end{figure}

Next, we demonstrate that the single-level algorithm \Cref{e:level} can be applied to $\tilde{A}$ recursively.
In order to factorize $\tilde{A}$ (approximately), we consider the remaining skeleton points at the coarse level $L-1$. Let every box $\B_i$ own the skeleton points of its children, i.e., $\cup_{k \in \C(\B_i)} \S_k$. See \Cref{f:merge} for a pictorial illustration. We obtain a partitioning of $\tilde{A}$ as we did for $A$ (a leaf box $\B_i$ owns points physically lying inside it).


In order to continue applying the strong skeletonization operator to every box at level $L-1$, we need to first verify \Cref{a:bf} still holds. For the skeletons $\S_i$ in box $\B_i$, the modified interactions (submatrices) only exist between $\S_i$ and $\cup_{k \in \M(\B_i)} \S_k$. At level $L-1$, it can be shown that the parent of $\B_i$ is the neighbor of the parent of any box in $\M(\B_i)$. As a result, \Cref{a:bf} holds, and thus we can repeat the previous approach to apply the strong skeletonization operator to every box at level $L-1$. By induction, we know that

\begin{theorem}
\Cref{a:bf} holds at all levels $\ell = L, L-1, \ldots, 2, 1$\footnote{There is no far-field for any box at level 1, so the theorem holds trivially.}.
\end{theorem}

We refer to~\cite[Section 3.3.2]{minden2017recursive} for a formal proof of the above theorem. The implication of the theorem is that \Cref{r:nn} and \Cref{r:t} both hold at all levels. Therefore, we obtain a multi-level algorithm:
\begin{equation} \label{e:alg}
Z(A;  \,
\underbrace{\B_1,\B_2,\ldots, \B_{n_L}}_{\text{level } L}, \,
\ldots, \,
\underbrace{\B_1,\B_2,\ldots, \B_{n_{1}}}_{\text{level } 1} \,
) \approx 
I
.
\end{equation}
As explained in \Cref{s:nn}, the inverse of a strong skeletonization operator can be applied efficiently (without forming the  inverse operator explicitly). We summarize the complete algorithm in \Cref{a:fact}.

As mentioned earlier, other low-rank approximation methods can also be used to construct similar factorizations~\cite{ambikasaran2014inverse,coulier2017inverse}. The same framework can be applied to solving large \emph{sparse} linear systems with minor modifications~\cite{pouransari2017fast,sushnikova2018compress}. While the described framework leverages the so-called strong admissibility (compression of far-field), a few methods based on the weak admissibility (compression of both far-field and near-field) can also construct efficient factorizations~\cite{corona2015n,ho2016hierarchical}.

\begin{algorithm}
\caption{\texttt{sequential factorization}}
\label{a:fact}
\begin{algorithmic}[1]
\Require{Kernel function $K$, points $\X = \{x_i\}_{i=1}^N \in \mathbb{R}^{2}$.}
\Ensure Factorization in \cref{e:alg}.
%
\State  Compute a hierarchical partition of $\X$. For every box $\B$, we obtain its
{neighbors $\N$ and distance-2 neighbors $\M$ (\cref{d:m}).} 
%
\For{$\ell = L, L-1, \ldots, 1$}  \Comment{\hfill \cyan{// bottom-up sweep}}
\For{$i = 1, 2, \ldots, n_\ell$} \Comment{\hfill \cyan{// a given ordering of boxes}}
%
%
\State \Call{Strong\_Skeletonization}{box $\B_i$ at level $\ell$}
\EndFor
\EndFor
\Statex
\end{algorithmic} 

\begin{algorithmic}[1]
\Function{Strong\_Skeletonization}{$\B$}
\State Read matrices $A_{\M, \B}$ and $A_{\B, \M}$, and compute the ID compression. 
%
%
\Comment{\hfill \cyan{// \cref{s:comp}}}
%
\State Read matrices $A_{\R, \R}$, $A_{\N, \R}$, and $A_{\R, \N}$; and update matrix $A_{\N, \N}$. 
%
\Comment{\hfill \cyan{//  \cref{s:nn}}}
\EndFunction

%

\end{algorithmic}
\end{algorithm}

\subsection{Solution/applying the inverse of factorization}

The solution phase follows the standard forward and backward substitution procedures besides applying the inverse of the sparsification operator \Cref{e:s}. The  algorithm  consists of an upward and a downward level-by-level traversal of the tree $\T$. During the upward pass, all boxes at the same level are visited according to the order of factorization. For every box $\B$, we first update a subset of the right-hand-side (RHS) $b(\B)$ and then update $b(\N)$ associated with neighbor boxes of $\B$ in our  in-place algorithm. The operations are ``reversed'' during the downward pass, where we read neighbor data $b(\N)$ and update local data $b(\B)$ for every box $\B$.

\section{Distributed-memory parallel algorithm}

\begin{figure}
\begin{center}
\includegraphics[width=\linewidth]{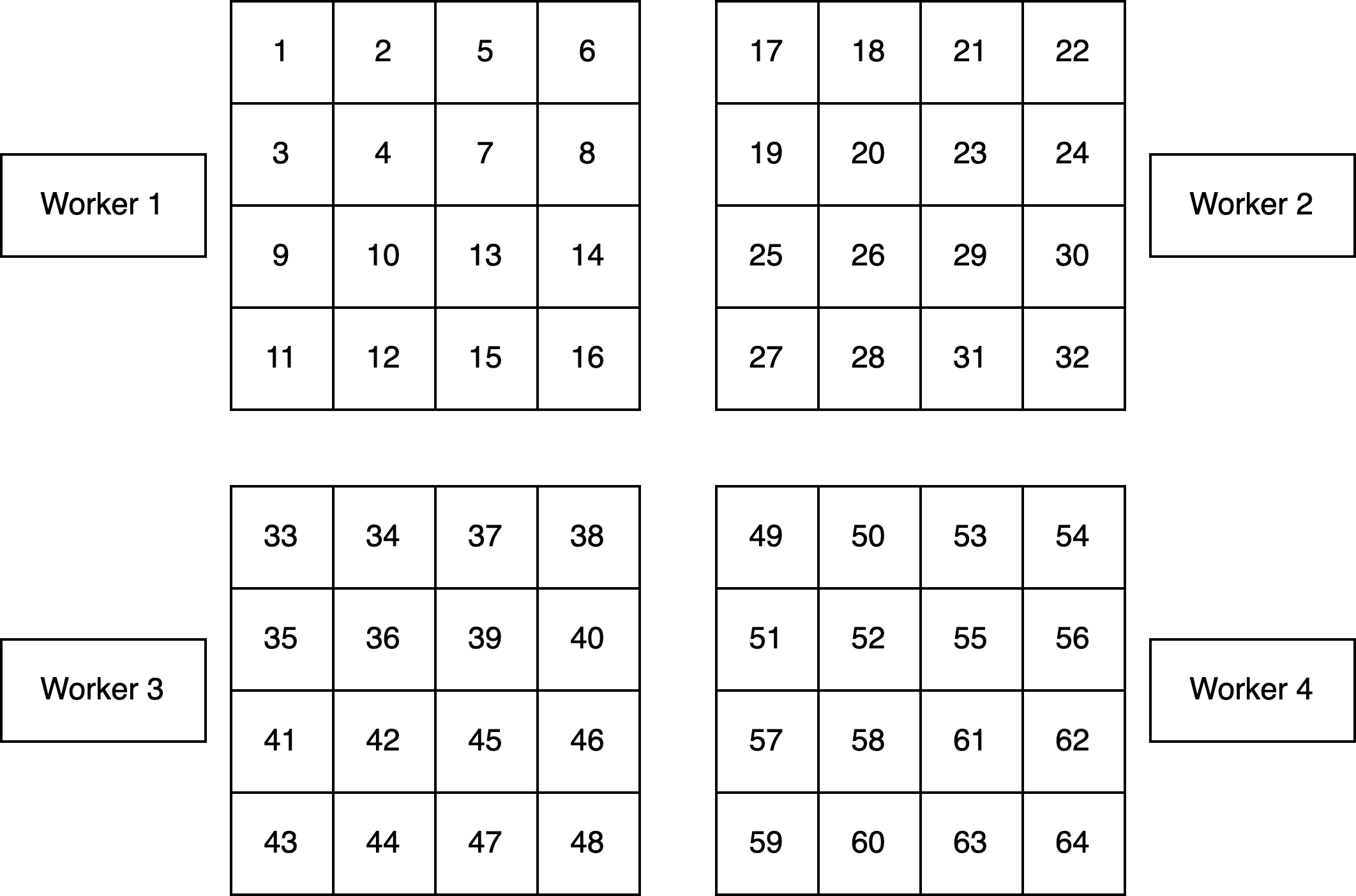}
\end{center}
\caption{Distribution of our data structure associated with boxes among all  worker processes. The data structure is stored using distributed arrays.
}
\label{f:partition}
\end{figure}

In the following, we analyze data dependency at the leaf level based on previous results summarized in \Cref{r:t,r:nn}, but the analysis also holds at other levels in the tree  (see \Cref{s:ml}). In \texttt{Julia}, a master process (process 1) provides the interactive prompt and coordinates worker processes to conduct parallel operations.  To store our data structures such as a list of required submatrices for every box, we used distributed arrays from \texttt{DistributedArrays.jl}.  A distributed array has the property that a process can make a fast local access but  has only read permission for a remote access. 
Our parallel algorithm starts by constructing  local data structures on each worker and assembling them into distributed arrays.

At the leaf level, all boxes are distributed evenly among all worker processes. See \Cref{f:partition} for an example. On every worker $\W_i$, the boxes are  classified into two groups: (1) ``interior boxes'', whose neighbors  are on the same worker, and (2) ``boundary boxes'', whose neighbors are on $\W_i$ and at least one other worker $\W_j$. In \Cref{f:partition}, boundary boxes on worker 1 include 6, 8, 14, 16, 11, 12, and 15; and the rest are all interior boxes. Notice that if the total number of boxes is large, the number of interior boxes dominates.

Consider a pair of boxes $\B_i$ and $\B_j$ at the leaf level. Suppose the two boxes are on two different workers $\W_i$ and $\W_j$, respectively. (Otherwise, the two boxes will be processed sequentially.) Consider the distance ${d}(\B_i, \B_j)$ between the two boxes,\footnote{${d}(\B_i, \B_j) = \max(|x_i-x_j|, |y_i-y_j|)/\ell$, where $(x_i, y_i)$ and $(x_j, y_j)$ are the centers of $\B_i$ and $\B_j$, respectively, and $\ell$ is the side length of  boxes at the same level as $\B_i$ and $\B_j$.}
where 
\begin{enumerate}
\item 
$\text{dist}(\B_i, \B_j) = 1$ if $\B_i$ and $\B_j$ are neighbors/adjacent (e.g., boxes 6 and 17 in \Cref{f:partition}). In fact, both boxes must be boundary boxes. We discuss the parallel algorithm for handling all boundary boxes in \Cref{sec: fac_scheduling}.

\item
$\text{dist}(\B_i, \B_j) = 2$ if $\B_i$ and $\B_j$ are not neighbors but share a common neighbor, i.e., they are distance-2 neighbors as defined in \Cref{d:m} (e.g., boxes 5 and 17 in \Cref{f:partition}). In fact, one of the two boxes must be a boundary box and the other one must be an interior box. If one of them, say $\B_i$, is processed first, it requires accessing the corresponding submatrices $A_{\B_i, \B_j}$ and $A_{\B_j, \B_i}$ (see \Cref{r:t}) but will not update submatrices in the rows/columns of $\B_j$  (see \Cref{r:nn}).

\item
$\text{dist}(\B_i, \B_j) > 2$ if $\B_i$ and $\B_j$ do not share any common neighbors (e.g., boxes 5 and 18 in \Cref{f:partition}). In this case, the two boxes can be processed in parallel according to \Cref{r:t,r:nn}.

\end{enumerate}

Our parallel algorithm for the described matrix factorization employs a level-by-level traversal of the quad-tree $\mathcal{T}$. At every level, we have three steps: (1) processing the massive interior boxes in parallel, (2) processing the remaining boundary boxes, and (3) creating the next/coaser level. 
The pseudocode is shown in \Cref{alg:FMM:Parallel}, and we explain implementation details in the following sections.
\begin{algorithm}
\caption{\texttt{parallel factorization}}
\label{alg:FMM:Parallel}
\begin{algorithmic}[1]

\For{$\text{tree level } \ell = L, L-1, \ldots, 1$} 
\State All workers factorize interior boxes at level $\ell$. 
\For {$\text{color } i = 1, 2, 3, 4$}  
\State Workers with color $i$ factorize boundary boxes at level $\ell$.
\State Workers with color $i$ send data to neighbors.
\EndFor
\State All workers construct  level  $\ell-1$.
\EndFor

\end{algorithmic}
\end{algorithm}

\subsection{Interior boxes}

For the first step in our parallel algorithm at every level, notice that if a pair of interior boxes $\B_i$ and $\B_j$ are on two different workers, then $\text{dist}(\B_i, \B_j) > 2$ (e.g., boxes 5 and 18 in \Cref{f:partition}). Therefore, all workers can apply the strong skeletonization operator (see \Cref{a:fact}) to factorize their respective interior boxes concurrently. In our implementation, every worker stores all necessary submatrices for this step. For example, in \Cref{f:partition} both worker 1 and worker 2 store copies of the submatrices $A_{5, 17}$ and $A_{17, 5}$. Recall  \Cref{r:t} that both box 5 and box 17 need the two submatrices for computing their IDs. 

{The distributed programming model in Julia (Distributed.jl) is similar to task-based programming in, e.g., OpenMP~\cite{ayguade2008design}, where the main process creates new processes (workers) and launches tasks to all workers through remote procedure calls. An example is the following:}
\begin{lstlisting}[language=Julia]
@sync for i in workerIDs
  @async remotecall_wait(factor, i, boxes)
end
\end{lstlisting}
where 
\texttt{remotecall\_wait} tells worker $i$ to factorize a set of boxes, the \texttt{async} macro means the task launch is asynchronous, and the \texttt{sync} macro on the loop makes sure that the master process waits for all workers to finish their tasks. (If \texttt{remotecall} is used instead, then the master process only waits for all workers to receive their tasks.)

\subsection{Boundary boxes}
\label{sec: fac_scheduling}

\begin{figure}
     \centering
     \begin{subfigure}[b]{0.25\textwidth}
         \centering
         \includegraphics[width=0.5\linewidth]{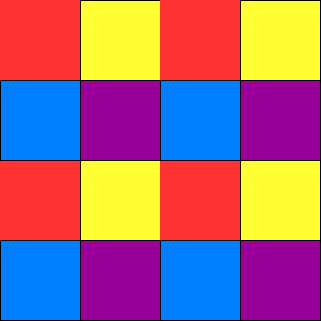}
     \end{subfigure}
%
     \caption{Coloring of a $4 \times 4$ process grid. For a 2D grid, we need 4 colors so that every process has a different color from its adjacent processes.
     }
     \label{f:color}
\end{figure}

After the first step, all workers are left with their boundary boxes to be processed. We assign colors to all workers such that any pair of adjacent processes have two different colors.  Although mature software packages exist for coloring an arbitrary  graph, we assume the boxes at every level form a grid (a full quad-tree $\T$), and thus a coloring of the process grid is straightforward as shown in \Cref{f:color}. Notice that 4 colors are needed for a 2D process grid regardless of the number of boxes on every worker.


Once the coloring is available, our parallel algorithm loops over all colors. For each color, the associated processors  need to fetch submatrices required for subsequent computation and then process their boundary boxes in parallel.
For the example in \Cref{f:color}, our algorithm will schedule 4 boxes, say with red color; wait for their completion; then schedule the next 4 boxes, say with yellow color; wait for their completion; and so on. Notice that if a pair of boundary boxes $\B_i$ and $\B_j$ are on two different processors with the same color, then $\text{dist}(\B_i, \B_j) > 2$ as long as every processor has at least $2\times 2 = 4$ boxes.

After all workers with a specific color $\alpha$ finish factorizing their boundary boxes, they need to send submatrices from Schur complement updates to their neighbors. 
In \texttt{Julia}, communication is ``one-sided,'' which means that only one process needs to be explicitly managed by the programmer. 
According to our best knowledge, remote direct memory access~\cite{https://doi.org/10.48550/arxiv.2208.01154} is not available in \texttt{Julia} for distributed computing. So we implemented explicit communication through remote procedural calls.
In particular, the master process launches \texttt{receive\_data} tasks on the neighbor workers as follows:
\begin{lstlisting}[language=Julia, mathescape=true]
@sync for i in neighborIDs
  @async 
  remotecall_wait(receive_data, i, $\alpha$)
end
\end{lstlisting}
Here, \texttt{receive\_data} is a function defined on all workers through the use of the \texttt{everywhere} macro as follows:
\begin{lstlisting}[language=Julia, mathescape=true]
@everywhere function receive_data($\alpha$)
  @sync for k in neighborWithColor[$\alpha$]
    @async 
    data=remotecall_wait(send_data, k, myID)
  end  
end
@everywhere function send_data(id)
  return UpdatesForWorker[id]
end
\end{lstlisting}
where the worker that needs data launches tasks on its neighbors who have color $\alpha$, and data will be returned in the form of a \texttt{future} object. Finally, the neighbor with color $\alpha$ returns/sends the requested data back, which is stored in a distributed array.

\subsection{Level transition}

Since the algorithm traverses  the quad-tree in a level-by-level fashion starting from the leaf level, it is important to keep track of the data structure when the algorithm transitions from one level to the next. The approach we take is to explicitly store the modified interactions for every box. As a result, whenever the algorithm changes level, it will need to reconstruct the data structure by regrouping interactions among accumulated points. An additional thing that we need to keep track of is where the coordinates of the given points are stored. At the start of the algorithm, the points are stored evenly distributed across all processes using a distributed array. However, as the algorithm progresses from one level to the next, two things will change. First, since each new level are constructed from the skeletonized points of the previous level, the points that are involved in the factorization of the new level are a subset of the previous points. Second, the number of processes involved in the new level may also decrease. Hence, this suggests that during the transition stage, we also need to reconstruct a distributed array by using only the relevant skeletonized points. 

\section{Complexity analysis} \label{s:analysis}

For simplicity, we make the following three assumptions. First, the discretization points $\{\bm{x}_i\}_{i=1}^N$ in \Cref{e:aij} are uniformly distributed in a square domain, so a hierarchical partitioning of the points correspond to a perfect quad-tree $\T$ (see \Cref{f:bnf}). Second, the numerical ranks of  the two submatrices $A_{\mathcal B, \mathcal F}$ and $A_{\mathcal F, \mathcal B}$ are constant for every cluster $\B$ (see \Cref{s:lr}), which is denoted by a scalar $r$. We further assume $r$ does not depend on the  problem size $N$ (see numerical results in \Cref{f:ranks}). Third,  the number of points per box at the leaf level is $\bigO(r)$, so the number of levels in the quad-tree $\T$ is $L = \bigO(\log(N / r)) = \bigO( \log N )$.

\subsection{Computational and memory complexity}

The costs of the strong RS algorithm were derived in~\cite[Section 3.3.4]{minden2017recursive}: the factorization cost $t_f = \bigO(N)$, the solve cost $t_s = \bigO(N)$, and the memory footprint also scales as $\bigO(N)$, where $N$ is the problem size.
%
Similar results for the inverse FMM algorithm can be found in~\cite[Section 2.5]{coulier2017inverse}. In the parallel algorithm, we uniformly partition the problem, the required work, and the memory footprint among all processors.

%
%

\subsection{Communication cost}

Notice the quad-tree $\T$ has $L = \bigO( \log N )$ levels, and all processors are organized as another quad-tree with $\bigO( \log p )$ levels. So every processor owns a subtree of $\bigO( \log (N/p) )$ levels.
In our parallel factorization algorithm, every processor sends a constant number of messages with a constant number of words for the first $\bigO( \log (N/p) )$ levels. Then, the remaining algorithm behaves as a parallel reduction among $p$ processors. As a result, the number of messages sent by every processor is $\bigO(\log N +  \log p )$, and the number of words moved is 
\begin{equation} \label{e:comm}
\bigO( \sqrt{{N}/{rp}} + \log p) = \bigO( \sqrt{N/p} + \log p)
\end{equation}
where $\sqrt{{N}/{rp}}$ is the number of boundary boxes on every processor.
For strong scaling experiments, i.e., $N = \bigO(1)$, the number of messages required is $\bigO(1)$, and the number of words is $\bigO(1/\sqrt{p} + \log p)$.
%
%
For weak scaling experiments, i.e., $N/p = \bigO(1)$, the number of messages required is $\bigO(\log p)$, and the number of words is $\bigO(\log p)$.
%


\section{Numerical results} \label{s:r}

In this section, we show benchmark results for solving two types of dense linear systems  $Ax=b$ that are associated with the free-space Green’s function for the Laplace equation in 2D and that for the Helmholtz equation in 2D, respectively. They represent kernel functions that are non-oscillatory and mildly oscillatory. 
For ease of setting up experiments, we discretized the integral equation \Cref{e:ie} on uniform grids  in the unit square. As a result, the matrix-vector product with dense matrix $A$ can be performed efficiently via the fast Fourier transform. (Otherwise, a fast summation algorithm such as the distributed-memory FMM is required.) \Cref{t:notations} summaries our notations for reporting numerical results.

\begin{table}
    \caption{\em  Notations used to report  results of solving $Ax=b$, where $A$ is from discretizing the integral equation \Cref{e:ie}, and $b$ is a standard uniform random vector.
    Timings are in seconds.
    }    
    \label{t:notations}
    \centering
    \begin{tabular}{ll} \toprule
    $N$ & size of matrix $A$, i.e., $A \in \mathbb{C}^{N \times N}$. \\
        $\varepsilon$ & tolerance for low-rank compression; see \Cref{definition:id}. Unless \\ & otherwise noted,   $\varepsilon = 10^{-6}$ is used. \\
    $p$ & number of processes. {In strong scaling tests ($N$ fixed),  we} \\ & 
    {started   with the minimum   number of compute nodes (processes)}  \\ & 
    {and increased the number of processes to at most 64 per node.} \\
    \midrule
    $t_{\text{fact}}$ & wall time of constructing the factorization \Cref{e:alg} in parallel. \\
    $t_{\text{solve}}$ & wall time of applying the inverse of the factorization in parallel. \\
    $t_{\text{comp}}$ & fraction of $t_{\text{fact}}$ or $t_{\text{solve}}$ spent on computation. \\
    $t_{\text{other}}$ & fraction of $t_{\text{fact}}$ or $t_{\text{solve}}$ spent on communication and overhead. \\ 
    \midrule
    relres & \textbf{rel}ative \textbf{res}idual, i.e., 
    $\|A \tilde{x} - b\|/\|b\|$, where $\tilde{x}$ is the output of \\ & our solver. \\
    $n_{it}$ & number of preconditioned CG or GMRES iterations to reach  \\ & $10^{-12}$ tolerance, where our solver is used as a preconditioner. \\
    \bottomrule
    \end{tabular}

\end{table}



All experiments were performed with \texttt{Julia}, version 1.9.4, on the CPU nodes of Perlmutter\footnote{\url{https://docs.nersc.gov/systems/perlmutter/architecture/}}, an HPE (Hewlett Packard Enterprise) Cray EX supercomputer. Each compute node has  two AMD EPYC 7763 (Milan) CPUs, 64 cores per CPU, and 512 GB of memory.
Since \texttt{Julia}  employs a just-in-time compiler that compiles a code before its first execution, we always ran a small problem size before timing our numerical experiments. 

\subsection{Laplace kernel}

Consider the following problem involving the free-space Green’s function for the 2D Laplace equation: solving the first-kind volume integral equation
\begin{equation} \label{e:int}
\int_{\Omega} K(\|\bm{x}-\bm{y}\|) u(\bm{y}) d \bm{y} = f(\bm{x}), \quad \bm{x} \in \Omega = [0,1]^2
\end{equation}
with the kernel function 
\begin{equation} \label{e:kernel_laplace}
K(\bm{x}_i, \bm{x}_j) = - \frac{1}{2\pi}\log\left(\|\bm{x}_i - \bm{x}_j \|\right). 
\end{equation}
We discretized \Cref{e:int} using piecewise-constant collocation on a $\sqrt{N} \times \sqrt{N}$ uniform grid $\X$. 
The resulting  linear system involves a dense matrix $A \in \mathbb{R}^{N\times N}$:  the off-diagonal entries are 
\begin{equation} \label{e:a1}
A_{i,j} = - \frac{h^2}{2\pi}\log\left(\|\bm{x}_i - \bm{x}_j \|\right), \quad \forall i\not=j,
\end{equation}
where $\bm{x}_i, \bm{x}_j \in \X$ and $h=1/\sqrt{N}$ is the grid spacing; and the  diagonal entries are given by
\begin{equation} \label{e:a1i}
A_{i,i} = 
 \int_{-h/2}^{h/2} \int_{-h/2}^{h/2} 
- \frac{1}{2\pi}\log\left( \| \bm{x} \| \right) dx_1 \, dx_2,  \quad \forall i,
\end{equation}
where $\bm{x}=(x_1,x_2) \in \mathbb{R}^2$. We evaluated \Cref{e:a1i}  numerically using an adaptive quadrature (\texttt{dblquad} from \texttt{MultiQuad.jl}).

\begin{figure}
     \centering
     \begin{subfigure}{0.24\textwidth}
         \centering
         \includegraphics[width=\textwidth]{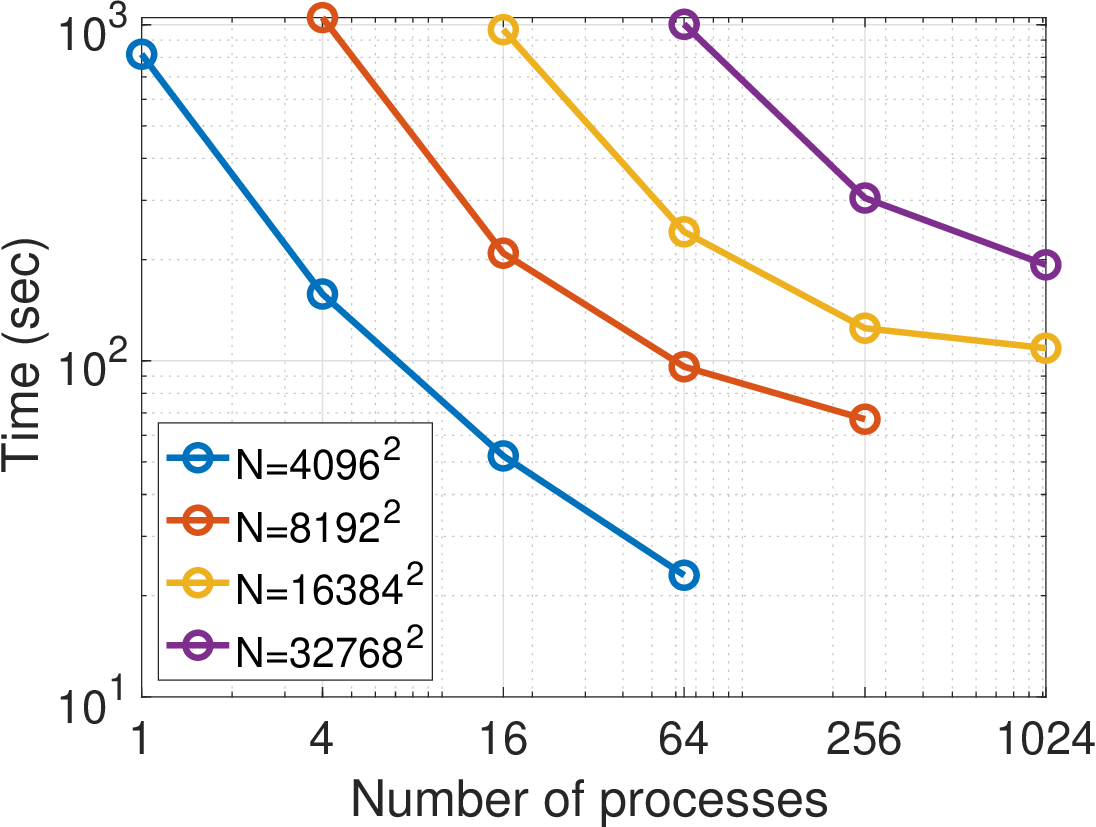}
         \caption{Strong scaling ($N$ fixed)}
     \end{subfigure}
	%
     \hfill
     \begin{subfigure}{0.24\textwidth}
         \centering
         \includegraphics[width=\textwidth]{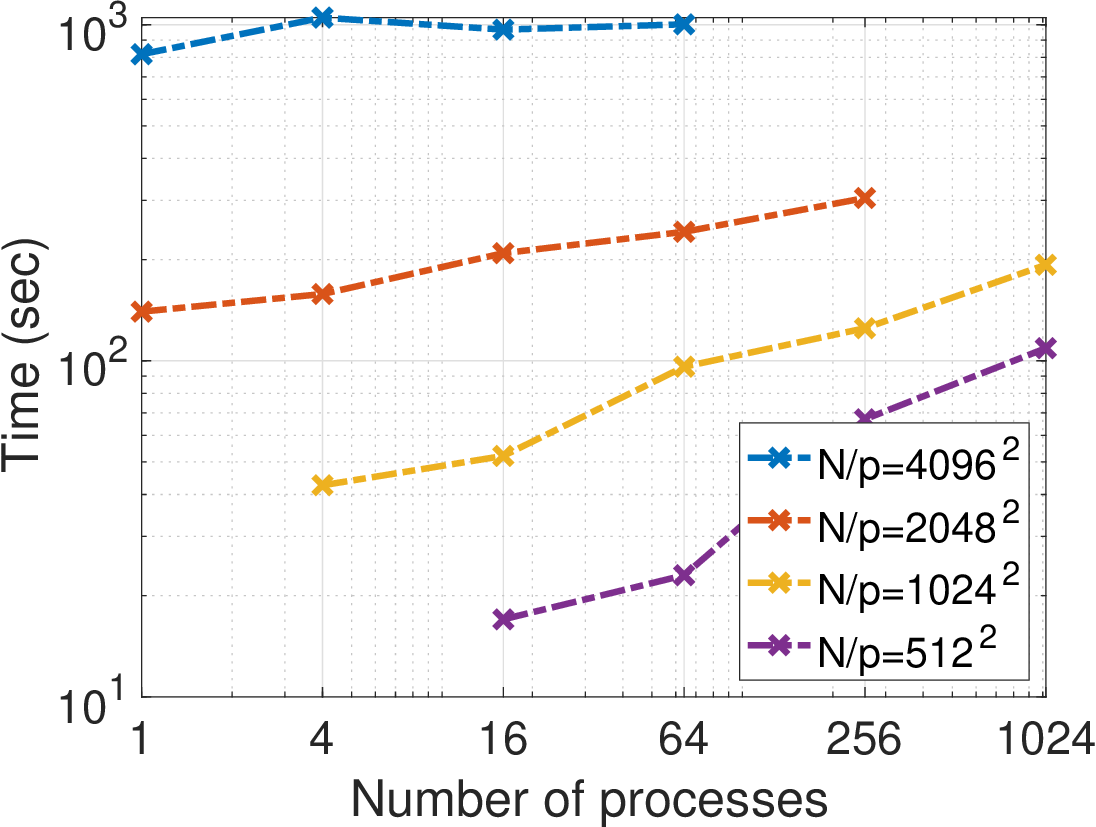}
         \caption{Weak scaling ($N/p$ fixed)}
     \end{subfigure}
     %
    \caption{{Scalability} of the factorization time $t_{\text{fact}}$   in \Cref{t:laplace}. 
    }
    \label{f:laplace}
\end{figure}

\Cref{t:laplace,t:laplace2} show the relevant numerical results, and \Cref{f:laplace} shows the parallel scalability of the factorization time. We make the following observations:
\begin{enumerate}
\item
We were able to solve a large problem size of $N=32\,768^2 \approx 1$ billion using 1024 processes in less than five minutes (200 seconds for factorization and around $11 \times 6 = 66$ seconds for six PCG iterations).

\item
The  approximate factorizations constructed in parallel required almost constant  number  of PCG iterations to arrive at a tolerance of $10^{-12}$. By contrast, the number of CG iterations is approximately $5 \sqrt{N}$ without any preconditioners. 

\item
The solve time for a single RHS is  much smaller than the factorization time. 



\end{enumerate}


\begin{table}
    \caption{\em Runtime for solving dense linear systems associated with the 2D Laplace kernel \Cref{e:kernel_laplace}.
    The tolerance for low-rank compression is $\varepsilon = 10^{-6}$, and the accuracy of our solver is shown in \Cref{t:laplace2}.
    }    
    \label{t:laplace}
    \centering
    \begin{tabular}{p{0.6cm}p{0.2cm}H|rrr|rrrH} \toprule
    \multirow{2}{*}{$N$}  &  \multirow{2}{*}{$p$}  &  \multirow{2}{*}{$M$}  
      & \multicolumn{3}{c|}{factorization time}  & \multicolumn{3}{c}{solve time ({one iteration})} & \multirow{2}{*}{$n_{\text{it}}$}  \\
      &  & & $t_{\text{fact}} = $ & $t_{\text{comp}} + $ & $t_{\text{other}}$ & $t_{\text{solve}} = $ & $t_{\text{comp}} + $ & $t_{\text{other}}$ &  \\ \midrule 
        %
$2048^2$ & 1 & 1 & 140 & 126 & 14 & 4.12 & 3.82 & 0.30 & 4 \\
& 4 & 1 & 42.6 & 36.9 & 5.7 & 1.50 & 1.25 & 0.25 & 4 \\ 
& 16 & 1 & 17.0 & 13.2 & 3.80 & 0.77 & 0.55 & 0.22 & 5 \\ 
& 64 & 1 & 10.6 & 7.6 & 3.0 & 0.64 & 0.37 & 0.27 & \\
\midrule
$4096^2$ & 1 & 1 & 817 & 719 & 98 & 36.36 & 32.30 & 4.06 & 5 \\
  & 4 & 1 & 158 & 141 & 17 & 6.19 & 5.47 & 0.72 & 5 \\
 & 16 & 1 & 52.1 & 42.2 & 9.9 & 2.06 & 1.63 & 0.43 & 5 \\ 
 & 64 & 1 & 23.0 & 17.1 & 5.9 & 1.41 & 1.00 & 0.41 & \\
 \midrule
$8192^2$ & 4 & 4 & 1050 & 847 & 203 & 54.11 & 48.25 & 5.86 & 5 \\
  & 16 & 4 & 209 & 147 & 62 & 7.19 & 6.12 & 1.07 & 5  \\
  & 64 & 4 & 96 & 46 & 50 & 3.69 & 2.81 & 0.88 & 5  \\
  & 256 & 4 & 67 & 21 & 46 & 2.96 & 1.63 & 1.33 & 5 \\ 
  \midrule
$16384^2$ & 16 & 16 & 967 & 749 & 218 & 30.63 & 27.20 & 3.43 & 5 \\
  & 64 & 16 & 242 & 158 & 84 & 10.99 & 8.18 & 2.81 & 5 \\
  & 256 & 16 & 125 & 51 & 74 & 5.39 & 4.03 & 1.36 & 5  \\
  & 1024 & 16 & 109 & 26 & 83 & 6.18 & 2.85 & 3.33 & 5 \\
    \midrule
{$32768^2$} 
& 64 & 64 & 1003 & 771 & 232 & 35.83 & 29.23 & 6.60 & \\
& 256 & 64 & 305 & 168 & 137 & 9.34 & 7.23 & 2.11 & \\
& 1024 & 64 & 193 & 59 & 134 & 11.29 & 5.81 & 5.48 & \\
    \bottomrule
    \end{tabular}
\end{table}

\begin{table}
    \caption{\em {Accuracies} for solving dense linear systems associated with the 2D Laplace kernel \Cref{e:kernel_laplace}.
    $\varepsilon$ denotes the tolerance for low-rank compression. 
    }    
    \label{t:laplace2}
    \centering
    \begin{tabular}{crrrrlc} \toprule
   $\epsilon$ & $N$ & $p$ & $t_{\text{fact}}$ & $t_{\text{solve}}$ & relres & $n_{\text{it}}$ 
   \\ \midrule
    %
    %
    \multirow{4}{*}{$10^{-6}$}
     & $2048^2$ & 1 & 139 & 4.12 & 1.11e-4 & 4 \\
     & $4096^2$ & 4 & 158 & 6.61 & 2.63e-4 & 5 \\
     & $8192^2$ & 16 & 209 & 9.85 & 4.48e-4 & 5  \\
     & $16384^2$ & 64 & 242 & 12.01 & 6.57e-4 & 6  \\
     \midrule
    \multirow{4}{*}{$10^{-9}$}
     & $2048^2$ & 1 & 225 & 4.69 & 1.31e-7 & 2 \\
     & $4096^2$ & 4 & 262 & 6.09 & 2.14e-7 & 2 \\
     & $8192^2$ & 16 & 374 & 10.73 & 2.48e-7 & 3  \\
     & $16384^2$ & 64 & 452 & 14.17 & 5.39e-7 & 3  \\
     \midrule
    \multirow{4}{*}{$10^{-12}$}
     & $2048^2$ & 1 & 446 & 5.85 & 1.44e-10 & 2 \\
     & $4096^2$ & 4 & 642 & 7.42 & 2.58e-10  & 2  \\
     & $8192^2$ & 16 & 867 & 9.25 & 4.47e-10 & 2   \\
     & $16384^2$ & 64 & 1048 & 16.45 & 5.58e-10 & 2  \\
    \bottomrule
    \end{tabular}
\end{table}

\subsection{Helmholtz Kernel}

Consider the next problem involving the free-space Green’s function for the 2D Helmholtz equation: solving the second-kind volume integral equation
\begin{equation} \label{e:int2}
    \sigma(\bm{x}) + \kappa^2 b(\bm{x}) \int_{\Omega} {K}(\|\bm{x} - \bm{y}\|) \sigma(\bm{y}) d \bm{y} =  - \kappa^2 b(\bm{x}) \, u_{\text{in}}(\bm{x})  
\end{equation}
where $\bm{x} \in \Omega = [0, 1]^2$,
known as the Lippmann-Schwinger equation, a reformulation of the variable coefficient Helmholtz equation (see, e.g., ~\cite[Section 11.2]{martinsson2019fast}) that models, e.g., acoustic wave propagation in a medium with a variable wave speed. Here $u_{\text{in}}(\bm{x})$ is the incoming wave with frequency $\kappa$; the ``scattering potential'' $0 < b(\bm{x}) \le 1$ is a known smooth function compactly supported inside $\Omega$; and the kernel function 
\begin{equation} \label{e:kernel_helmholtz}
K(\bm{x}_i, \bm{x}_j) = \frac{i}{4}H_{0}^{(1)}(\kappa \|\bm{x}_i - \bm{x}_j \|),
\end{equation}
where $H_0^{(1)}$ is the zero-order first-kind Hankel function.

We symmetrized \Cref{e:int2}  via change of variable $\mu(\bm{x}) = \sigma(\bm{x}) / \sqrt{b(\bm{x})}$   and  
applied piecewise-constant collocation on a $\sqrt{N} \times \sqrt{N}$ uniform grid $\X$.   The resulting \emph{indefinite complex} linear system involves a dense matrix $A \in \mathbb{C}^{N\times N}$: the off-diagonal entries are 
\begin{equation} \label{e:b1}
A_{i,j} = h^2 \kappa^2 \sqrt{b(\bm{x}_i) b(\bm{x}_j)} \left( \frac{i}{4}H_{0}^{(1)}(\kappa \|\bm{x}_i - \bm{x}_j\|) \right)
\end{equation}
where $\bm{x}_i, \bm{x}_j \in \X$, and $h=1/\sqrt{N}$ is the grid spacing; and the  diagonal entries are given by
\begin{equation} \label{e:b2}
A_{i,i} = 1 +  \kappa^2 b(\bm{x}_i)
 \int_{-h/2}^{h/2} \int_{-h/2}^{h/2} 
\frac{i}{4}H_{0}^{(1)} \left( \kappa \| \bm{x} \| \right)
dx_1 \, dx_2  
\end{equation}
where $\bm{x}=(x_1,x_2) \in \mathbb{R}^2$. We evaluated  \Cref{e:b2} numerically using an adaptive quadrature (\texttt{dblquad} from \texttt{MultiQuad.jl}).

\begin{figure}
\centering
     \begin{subfigure}{0.24\textwidth}
         \centering
         \includegraphics[width=\textwidth]{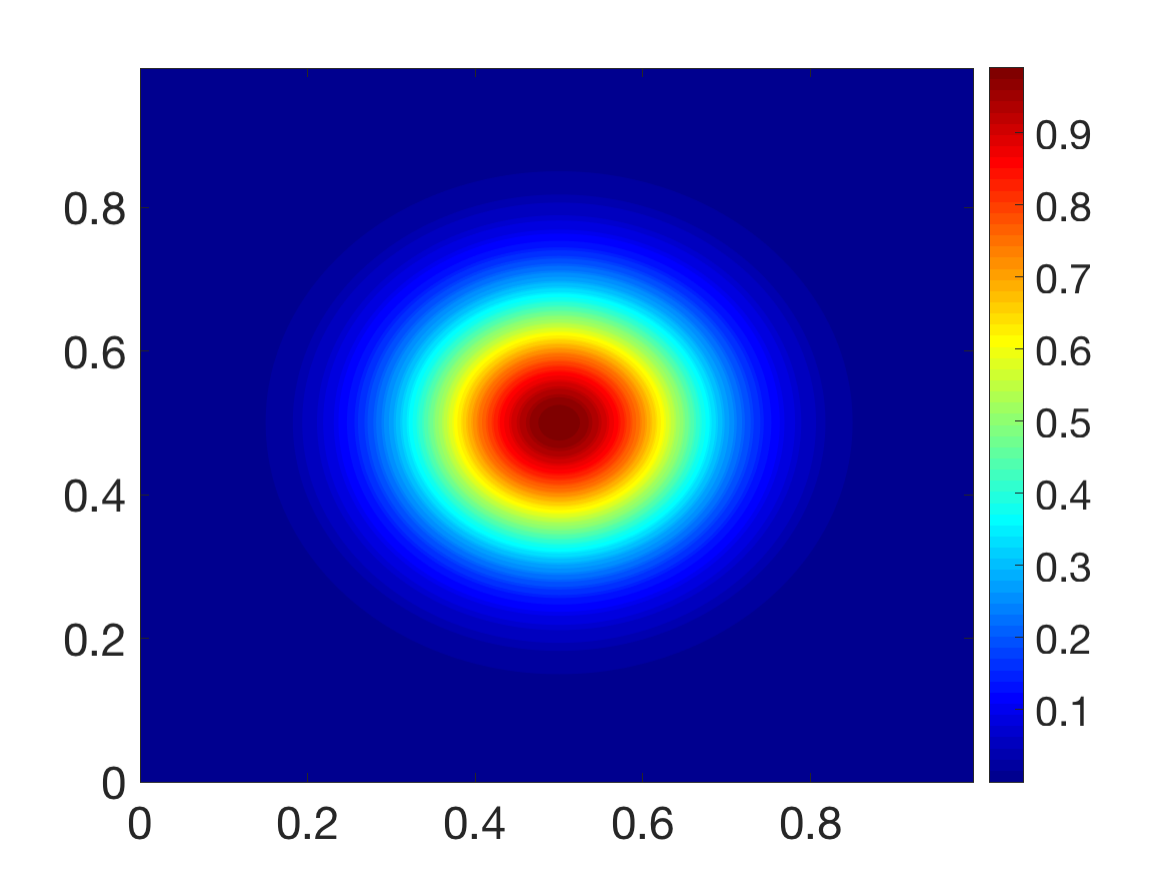}
         \caption{Scattering potential}
         \label{f:potential}
     \end{subfigure}
     \hfill
     \begin{subfigure}{0.24\textwidth}
         \centering
         \includegraphics[width=\textwidth]{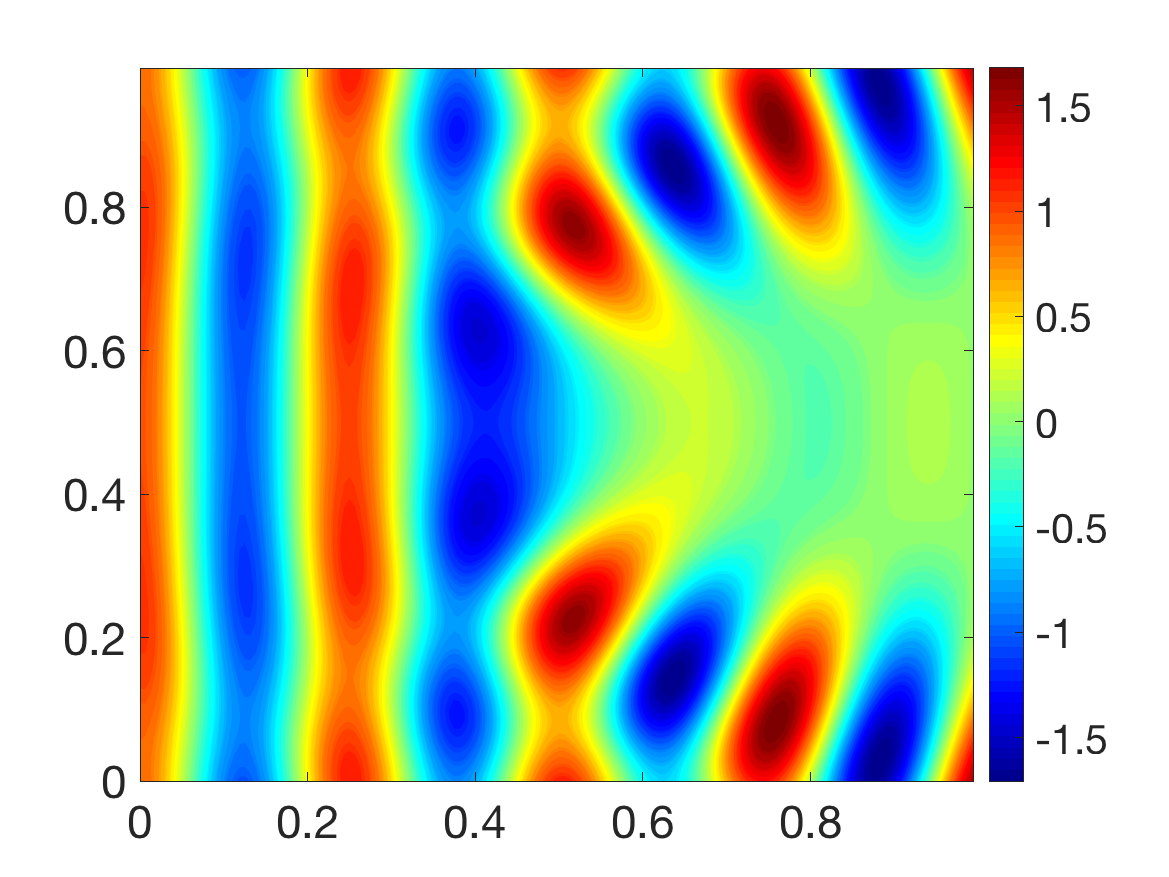}
         \caption{Total field}
         \label{f:scattering}
     \end{subfigure}
     %
     %
\caption{Gaussian bump scattering potential and the total field for an incoming plan wave  traveling from left to  right.
}
\label{f:helmholtz}
\end{figure}

For the following experiments, we use a Gaussian bump scattering potential $b(\bm{x}) = e^{-32 \|\bm{x} - \bm{c}\|^2}$ with the center $\bm{c} = (\frac{1}{2}, \frac{1}{2})$, as shown in \Cref{f:potential}. For an incoming plan wave $u_{\text{in}}(\bm{x})$ pointing to the right, the total field 
$
u(\bm{x}) = u_{\text{in}}(\bm{x}) + \int_{\Omega}   {K}(\|\bm{x}-\bm{y}\|) \sigma(\bm{y}) d \bm{y} 
$
is shown in \Cref{f:scattering}.

\begin{table}
    \caption{\em  Results for solving dense linear systems associated with the 2D Helmholtz kernel \Cref{e:kernel_helmholtz}  with \emph{fixed} frequency $\kappa = 25$.
    The tolerance for low-rank compression is $\varepsilon = 10^{-6}$, and the accuracy of our solver is shown in \Cref{t:cpp}. 
    }    
    \label{t:helmholtz}
    \centering
    \begin{tabular}{p{0.5cm} >{\raggedleft}p{0.2cm} H|rrr|rrrH} \toprule
    \multirow{2}{*}{$N$}  &  \multirow{2}{*}{$p$}  &  \multirow{2}{*}{$M$}  
      & \multicolumn{3}{c|}{factorization time}  & \multicolumn{3}{c}{solve time ({one iteration})} & \multirow{2}{*}{$n_{\text{it}}$}  \\
      &  & & $t_{\text{fact}} = $ & $t_{\text{comp}} + $ & $t_{\text{other}}$ & $t_{\text{solve}} = $ & $t_{\text{comp}} + $ & $t_{\text{other}}$ &  \\ \midrule 
        %
$1024^2$ & 1 & 1 & 315 & 300 & 15 & 1.77 & 1.60 & 0.17 & 3\\
   & 4 & 1 & 82 & 76 & 6 & 0.64 & 0.50 & 0.14 & 3\\
 & 16 & 1 & 27 & 23 & 4 & 0.47 & 0.24 & 0.23 & 3\\ 
 & 64 & 1 & 11 & 8 & 3 & 0.56 & 0.20 & 0.36 & 3\\
 \midrule
$2048^2$ & 1 & 1 & 1273 & 1212 & 61 & 6.74 & 6.43 & 0.31 & 3\\
   & 4 & 1 & 313 & 294 & 19 & 2.00 & 1.74 & 0.26 & 3\\
 & 16 & 1 & 98 & 88 & 10 & 1.08 & 0.85 & 0.23 & 3\\ 
 & 64 & 1 & 44 & 38 & 6 & 0.89 & 0.54 & 0.35 & 3\\
 \midrule
$4096^2$ & 1 & 1 & 5116 & 4873 & 243 & 30.69 & 29.37 & 1.32 \\
& 4 & 4 & 1237 & 1174 & 63 & 8.05  & 7.31 & 0.74 & 3\\
  & 16 & 16 & 340 & 314 & 26 & 3.59 & 2.79 & 0.80 & 3 \\
  & 64 & 16 & 117 & 103 & 14 & 2.02 & 1.60 & 0.42 & 3 \\
  & 256 & 64 & 81 & 44 & 37 & 1.78 & 0.97 & 0.81 & 3\\ 
  \midrule
$8192^2$ & 4 & 4 & 4958 & 4636 & 322 & 36.08 & 33.26 & 2.82 & 3 \\
  & 16 & 4 & 1292 & 1178 & 114 & 12.65 & 9.88 & 2.77 & 3\\
  & 64 & 4 & 394 & 318 & 76 & 5.28 & 3.85 & 1.43 & 3\\
  & 256 & 4 & 178 & 109 & 69 & 4.62 & 3.21 & 1.41 & 3\\
  & 1024 & 16 & 104 & 50 & 54 & 4.47 & 2.32 & 2.15 \\
  \midrule
{$16384^2$}    & 64 & 64 & 1369 & 1187 & 182 & 16.07 & 13.35 & 2.72 & 3 \\
  & 256 & 64 & 419 & 315 & 104 & 7.72 & 5.36 & 2.36 & 3\\
  & 1024 & 64 & 213 & 110 & 103 & 7.57 & 3.78 & 3.79 \\
  \midrule
{$32768^2$}  
  & 256 & 256 & 1420 & 1206 & 215 & 22.44 & 15.33 & 7.12 \\  
  & 1024 & 256 & 468 & 329 & 139 & 17.13 & 8.69 & 8.44 \\
  & 4096 & 256 & 289 & 131 & 158 & 28.28 & 9.88 & 18.40 \\
    \bottomrule
    \end{tabular}

\end{table}

\begin{figure}
\centering
     \begin{subfigure}{0.24\textwidth}
         \centering
         \includegraphics[width=\textwidth]{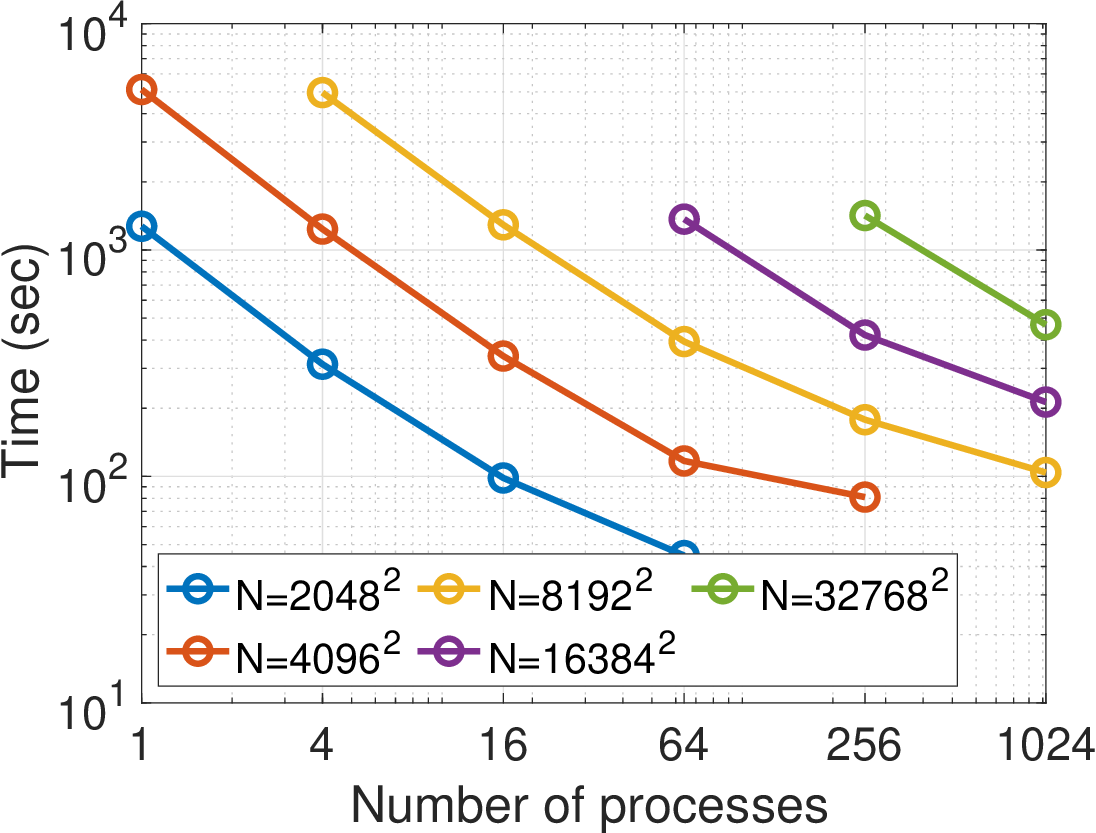}
         \caption{Strong scaling ($N$ fixed)}
     \end{subfigure}
     \hfill
     \begin{subfigure}{0.24\textwidth}
         \centering
         \includegraphics[width=\textwidth]{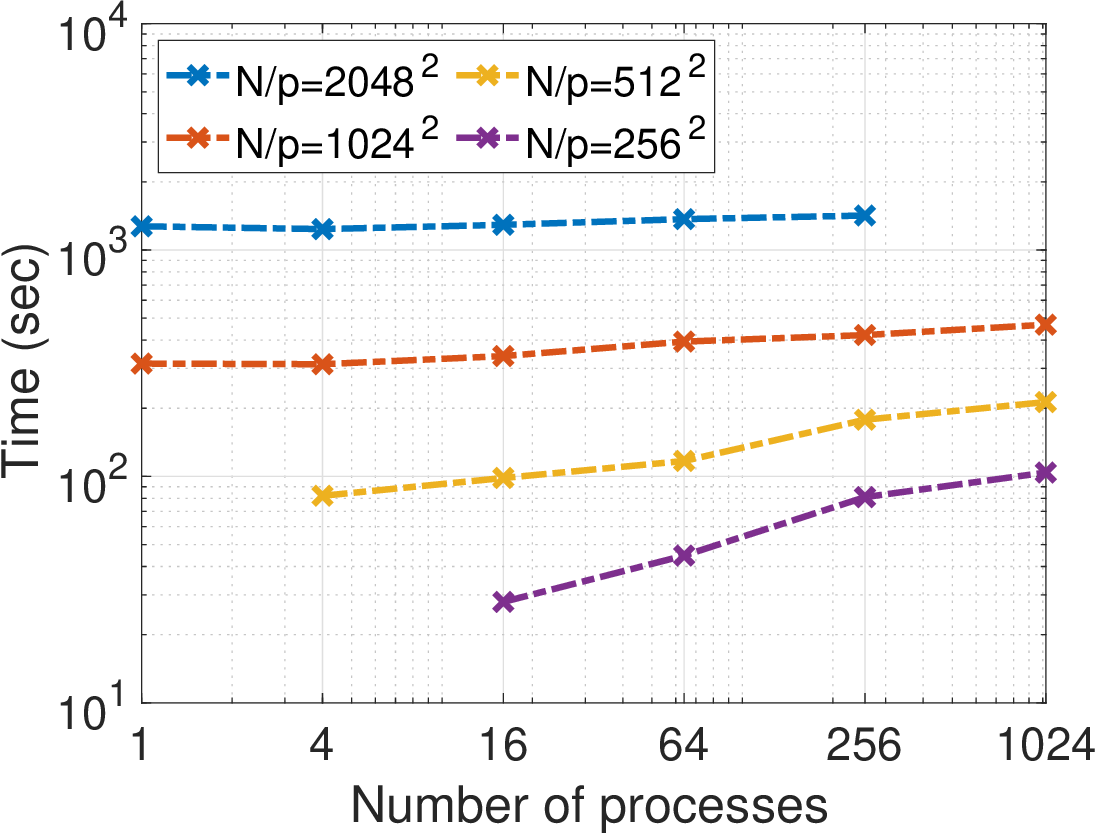}
         \caption{Weak scaling ($N/p$ fixed)}
     \end{subfigure}
\caption{{Scalability} of the factorization time $t_{\text{fact}}$ in \Cref{t:helmholtz}. 
}
\label{f:helmholtz}
\end{figure}


\subsubsection{Fixed frequency ($\kappa = 25$)}

\Cref{t:helmholtz} shows the relevant numerical results for a fixed frequency $\kappa = 25$, and \Cref{f:helmholtz} shows the parallel scalability of the factorization time. 
\begin{itemize}

\item
The factorization time is much longer than that for the Laplace kernel. The reason  is that an evaluation of {the (complex) Helmholtz kernel in \Cref{e:kernel_helmholtz} takes longer}.

\item
The parallel factorization algorithm achieves greater speedups compared to those for the Laplace kernel. The advantage of constructing a highly accurate approximation is clear: the solve time is much faster, suitable for situations where multiple RHSs need to be addressed.

\item 
Again, the  approximate factorizations constructed in parallel required a consistent number  (three) of GMRES iterations to arrive at a tolerance of $10^{-12}$.

\end{itemize}

\subsubsection{Increasing frequency ($\kappa = \pi \sqrt{N} /16 $, i.e., 32 points per wavelength)}

\begin{figure}
         \centering
         \includegraphics[trim=0 40 0 40,clip,width=0.49\textwidth]{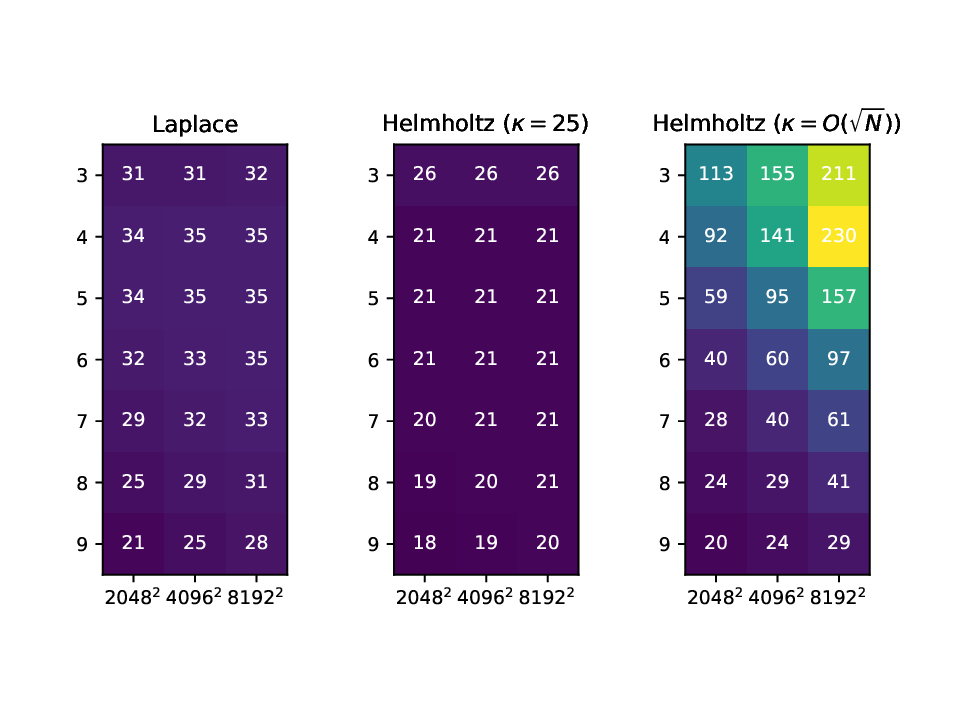}
    \caption{Numerical ranks of the 2D Laplace kernel \Cref{e:kernel_laplace} and the 2D Helmholtz kernel \Cref{e:kernel_helmholtz}. The y-axis stands for the tree level (see \Cref{f:tree}), where low-rank compression occurs starting from level 3. Every entry stands for the average rank of all boxes at a certain level.}
    \label{f:ranks}
\end{figure}

\Cref{t:h2} shows the relevant numerical results for incoming waves with increasing frequencies. 
\begin{itemize}

\item
The factorization time is increasingly longer than that in \Cref{t:helmholtz}, where the frequency was fixed. Recall that the numerical rank for the Helmholtz kernel is $\bigO(\kappa D)$; see \Cref{f:ranks}. As a result, the factorization requires $\bigO(\kappa^3) = \bigO(N^{3/2})$ operations. 

\item
As $\kappa$ increases in \Cref{e:int2}, the integral equation becomes increasingly more ill-conditioned. So is the discretized linear system, which requires more and more preconditioned GMRES iterations to converge. However, the savings from employing the preconditioner is clear: the number of GMRES iterations without any preconditioner is orders-of-magnitudes larger and grows rapidly.

\end{itemize}


\begin{table}
    \caption{\em 
    %
    Results for solving dense linear systems associated with the 2D Helmholtz kernel \Cref{e:kernel_helmholtz} with \emph{increasing} frequencies ($\kappa = \pi \sqrt{N} /16 $, i.e., 32 points per wavelength). 
    The last column $\tilde{n}_{\text{it}}$ shows the number of GMRES iterations (restart = 20) without any preconditioners. 
    }    \label{t:h2}
    \centering
    \begin{tabular}{ccHcrrrH|r} \toprule
    $N$ & $p$ & $M$ & ${\kappa}/({2\pi})$ & $t_{\text{fact}}$ & $t_{\text{solve}}$ & $n_{\text{it}}$ & ${||b - A\hat{x}||_2}/{||b||_2}$ & $\tilde{n}_{\text{it}}$ \\ \midrule
$1024^2$ & 1 & 1 & 32 & 365 & 1.83 &  3 & $3.02 \times 10^{-12}$ & 452 \\ 
$2048^2$ & 4 & 1 & 64 & 495 & 14.56 &  4 & $5.04 \times 10^{-12}$ & 1\,752 \\
$4096^2$ & 16 & 4 & 128 & 684 & 22.26 &  5 & $2.02 \times 10^{-10}$ & 8\,198 \\ 
$8192^2$ & 64 & 16 & 256 & 1138  & 29.74 & 12 & $2.26 \times 10^{-8}$ & $> 10\,000$ \\ 
    \bottomrule
    \end{tabular}
\end{table}

\subsection{{Comparison and 1-process-per-node run}}

\begin{figure}
\centering
     \begin{subfigure}{0.24\textwidth}
         \centering
         \includegraphics[width=\textwidth]{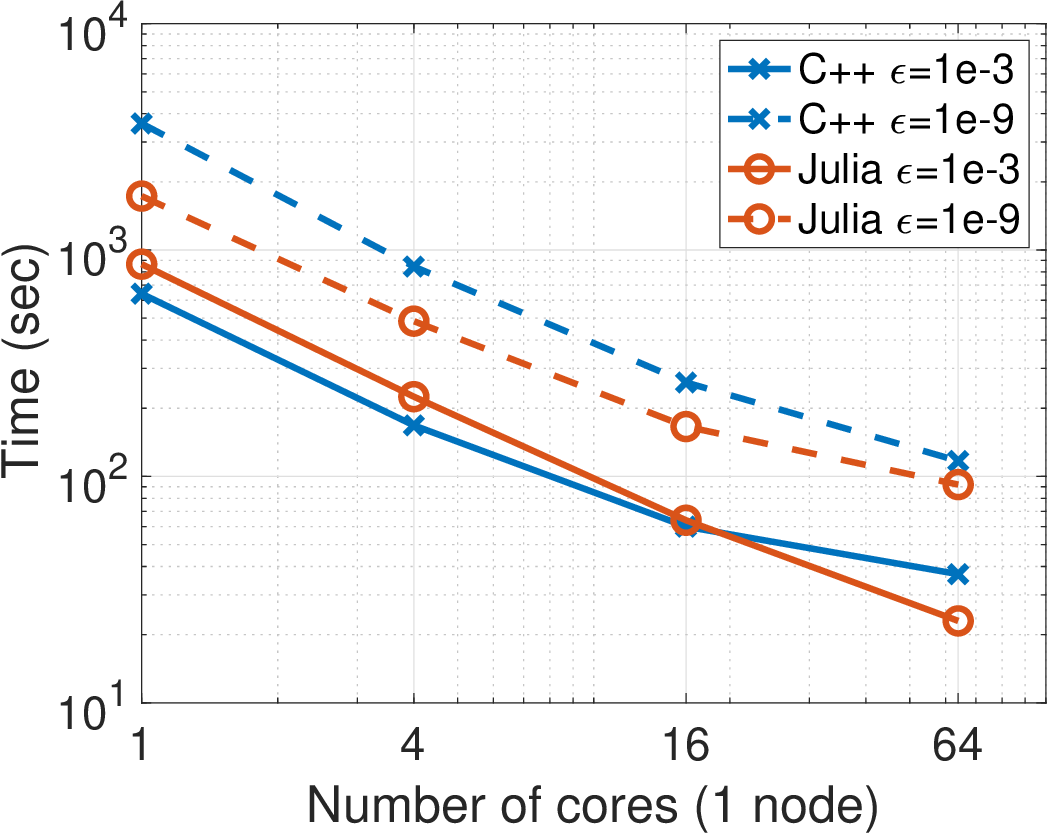}
     \end{subfigure}
     \hfill
     \begin{subfigure}{0.24\textwidth}
         \centering
         \includegraphics[width=\textwidth]{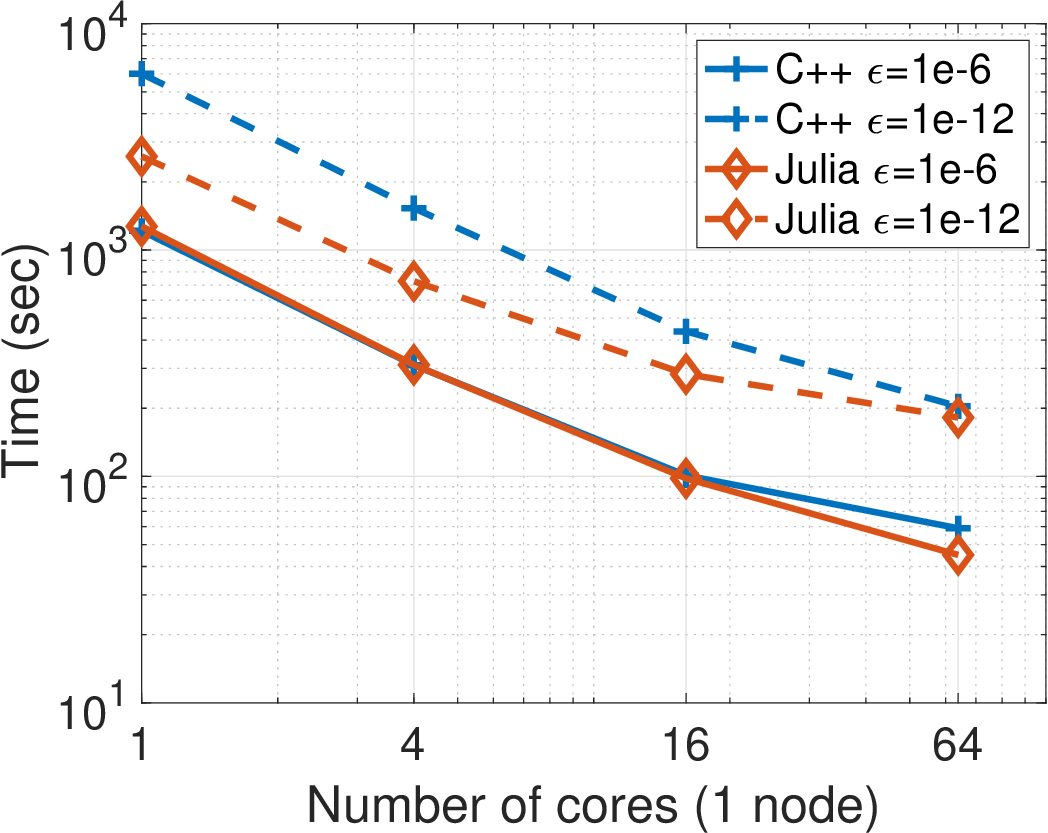}
     \end{subfigure}
     %
     %
\caption{
{Comparison of the factorization time ($t_{\text{fact}}$ in \Cref{t:cpp})  between a shared-memory solver in \texttt{C++} with \texttt{OpenMP}} (blue) and our distributed-memory solver in \texttt{Julia} (red). Both codes ran on one compute node with 1 process or 1 thread per core. 
    }
\label{f:cpp}
\end{figure}

We present two experiments on solving dense linear systems $Ax=b$ associated with the 2D Helmholtz kernel \Cref{e:kernel_helmholtz} with \emph{fixed} frequency $\kappa = 25$, where the matrix $A$ is defined in \Cref{e:b1,e:b2}.
The first one is a comparison between our distributed-memory solver implemented in \texttt{Julia} and a shared-memory  solver implemented in \texttt{C++} with \texttt{OpenMP}~\cite{ayguade2008design}. Both parallel solvers are based on  (sequential) RS-S~\cite{minden2017recursive}. The shared-memory solver follows ideas from Takahashi et al.~\cite{takahashi2020parallelization} and the parallel strategy briefly mentioned at the end of~\cite{minden2017recursive}: it colors all boxes at the same level, and every pair of neighbors has different colors.  By contrast, our distributed-memory solver colors only the boundary boxes on every process, and the boundary boxes on each process have the same color. The results of comparison are shown in \Cref{t:cpp,f:cpp}, where we observe that our \texttt{Julia} code performed similarly with the reference \texttt{C++} code when running on 64 cores of one compute node. The factorization of the Julia solver is faster but the solve time is slower with slightly worse accuracy.

\begin{table}
    \caption{\em  {Comparison between a shared-memory solver in \texttt{C++} with \texttt{OpenMP}} and our distributed-memory solver in \texttt{Julia}. Both codes ran on one compute node with one process or one thread per core. The matrix size $N=2048^2$.
    }    
    \label{t:cpp}
    \centering
    \begin{tabular}{H p{0.5cm} >{\raggedleft} p{0.15cm} | rrl | rrll} \toprule
    \multirow{2}{*}{$N$}  &  \multirow{2}{*}{$\epsilon$}  &  \multirow{2}{*}{$p$}  
      & \multicolumn{3}{c|}{C++ (reference)}  & \multicolumn{4}{c}{Julia (this paper)}   \\
      &  & & $t_{\text{fact}}$ & $t_{\text{solve}}$ & relres & $t_{\text{fact}} $ & $t_{\text{solve}}$ & relres & $n_{it}$  \\ \midrule 
      \multirow{16}{*}{$2048^2$}
       & \multirow{4}{*}{$10^{-3}$} & 1 & 641 & 5.47 & \multirow{4}{*}{8.5e-4}  & 866 & 5.63 & 1.9e-3 & 5\\
      & & 4 & 168 & 1.50 && 225 & 1.92 & 1.8e-3 & 6\\
      & & 16 & 60 & 0.46 && 64 & 0.76 & 1.8e-3 & 5\\
      & & 64 & 37 & 0.25 && 23 & 0.58 & 2.1e-3 & 5\\
      \midrule
       & \multirow{4}{*}{$10^{-6}$} & 1 & 1215 & 8.11 & \multirow{4}{*}{4.1e-7} & 1272 & 6.74 & 3.2e-6 & 3\\
      & & 4 & 310 & 2.17 && 313 & 2.00 & 4.2e-6 & 3\\
      & & 16 & 101 & 0.82 && 98 & 1.08 & 2.5e-6 & 3\\
      & & 64 & 59 & 0.39 && 45 & 0.89 & 3.5e-6 & 3\\
      \midrule
       & \multirow{4}{*}{$10^{-9}$} & 1 & 3621 & 13.18 & \multirow{4}{*}{9.1e-10} & 1730 & 7.11 & 3.3e-9 & 2 \\
      & & 4 & 845 & 3.22 && 485 & 3.06 & 1.3e-9 & 2\\
      & & 16 & 260 & 1.06 && 166  & 1.24  & 1.5e-9 & 2\\
      & & 64 & 117 & 0.76 && 92 & 1.03 & 2.1e-9 & 2\\    
      \midrule
       & \multirow{4}{*}{$10^{-12}$} & 1 & 6010 & 18.22 & \multirow{4}{*}{7.9e-13} & 2591 & 8.33 & 4.2e-12 & 2 \\
      & & 4 & 1531 & 5.18 && 728 & 3.59 & 6.3e-12 & 2\\
      & & 16 & 437 & 1.54 && 283 & 1.34 & 6.1e-12 & 2\\
      & & 64 & 204 & 1.02 && 182 & 1.30 & 8.7e-12 & 2\\      
    \bottomrule
    \end{tabular}

\end{table}

The second experiment investigates the extra costs for launching one process per compute node rather than launching multiple processes per compute node in previous sections. In particular, we reran a subset of experiments in \Cref{t:helmholtz} using as many compute nodes as the number of processes. Intuitively, we expect the new experiments to take longer due to more  communication through a network.  The results in \Cref{t:nodes} show that the extra wall-clock time is, however, negligible for these large-scale problems on  the Perlmutter supercomputer. 

\begin{table}
    \caption{\em Results for launching {1 process per compute node} for a subset of experiments in \Cref{t:helmholtz}. Here, the number of processes $p$ equals the number of compute nodes.}    
    \label{t:nodes}
    \centering
    \begin{tabular}{p{0.5cm} >{\raggedleft} p{0.2cm} H|rrr|rrrH} \toprule
    \multirow{2}{*}{$N$}  &  \multirow{2}{*}{nodes}  &  \multirow{2}{*}{$M$}  
      & \multicolumn{3}{c|}{factorization time}  & \multicolumn{3}{c}{solve time ({one iteration})} & \multirow{2}{*}{$n_{\text{it}}$}  \\
      &  & & $t_{\text{fact}} = $ & $t_{\text{comp}} + $ & $t_{\text{other}}$ & $t_{\text{solve}} = $ & $t_{\text{comp}} + $ & $t_{\text{other}}$ &  \\ \midrule 
%
$1024^2$   & 1 & 1 & 315 & 300 & 15 & 1.77 & 1.60 & 0.17 & 3\\
 \midrule
$2048^2$ & 1 & 1 & 1273 & 1212 & 61 & 6.74 & 6.43 & 0.31 & 3\\
   & 4 & 4 & 332 & 292 & 40 & 2.09 & 1.74 & 0.35 & 3\\
 \midrule
$4096^2$ & 1 & 1 & 5116 & 4873 & 243 & 30.69 & 29.37 & 1.32 \\ 
 & 4 & 4 & 1263 & 1160 & 103 & 8.86  & 7.95 & 0.91 & 3\\
  & 16 & 16 & 381 & 311 & 70 & 3.24 & 2.68 & 0.56 & 3 \\
  \midrule
$8192^2$ & 4 & 4 & 4958 & 4636 & 322 & 36.08 & 33.26 & 2.82 & 3 \\
  & 16 & 16 & 1333 & 1172 & 161 & 11.97 & 10.69 & 1.28 & 3\\
  & 64 & 64 & 397 & 315 & 82 & 4.38 & 3.35 & 1.03 & 3\\
  \midrule
$16384^2$  
  & 64 & 64 & 1369 & 1187 & 182 & 16.07 & 13.35 & 2.72 & 3\\
  & 256 & 256 & 421 & 314 & 107 & 6.97 & 4.67 & 2.30 & 3\\
  \midrule
  $32768^2$  
  & 256 & 256 & 1420 & 1206 & 215 & 22.45 & 15.33 & 7.12 \\
    \bottomrule
    \end{tabular}

\end{table}

\begin{figure}[H]
         \centering
         \includegraphics[width=0.25\textwidth]{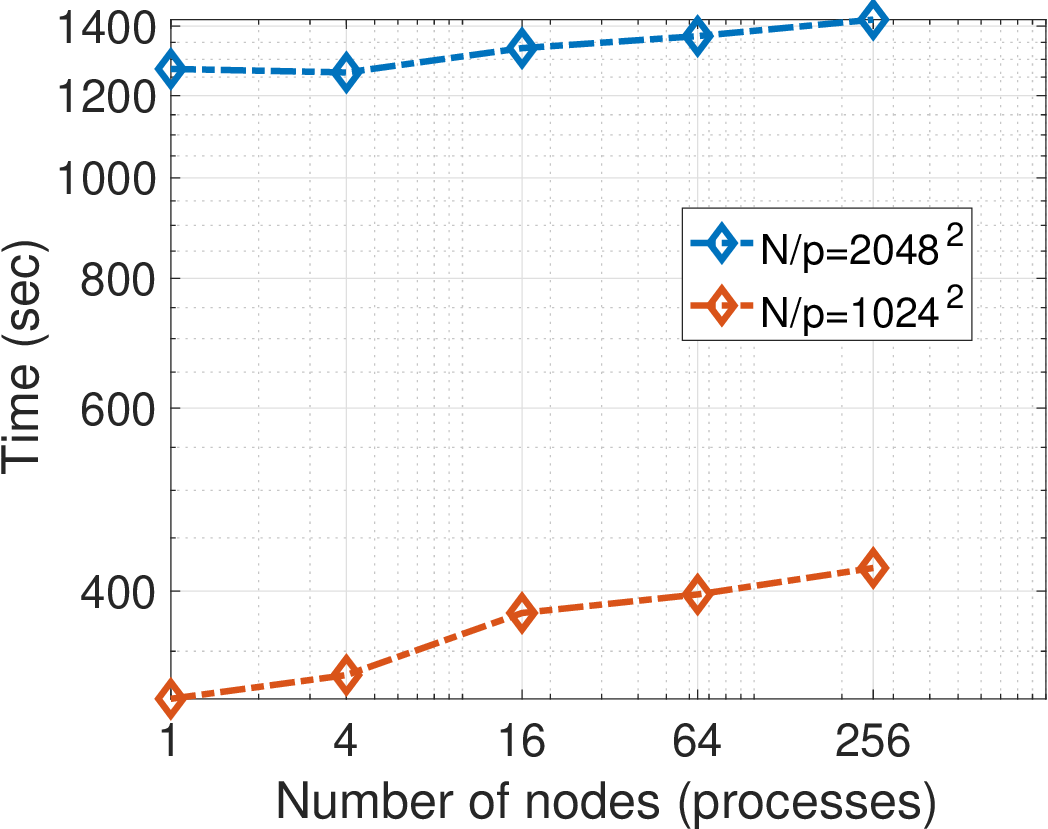}
         \caption{Weak scaling of the factorization time $t_{\text{fact}}$ in \Cref{t:nodes}.  {One process was launched on every compute node}.}
\label{f:scaling}
\end{figure}

\section{Conclusion}

We introduce a distributed-memory direct solver with $\bigO(N)$ complexity for solving dense linear systems arising from the discretization of  planar integral equations. Our algorithm follows from a data-dependency analysis of the RS-S method~\cite{minden2017recursive}. Our solver is implemented in \texttt{Julia}, and   numerical experiments are presented to show its parallel scalability.

\bibliographystyle{IEEEtran}
\bibliography{biblio}
\end{document}